\documentclass[reqno]{amsart}
 \usepackage{amsbsy,amssymb,amscd,amsfonts,latexsym,amstext,delarray,
 amsmath,graphicx, color}
 \usepackage[dvips]{epsfig}
\usepackage{hyperref}
\input xypic

\definecolor{Green}{rgb}{0.0,0.40,0.0}

\newtheorem{thm}{Theorem}[section]
\newtheorem{prop}[thm]{Proposition}
\newtheorem{cor}[thm]{Corollary}
\newtheorem{lem}[thm]{Lemma}

\newtheorem{defn}[thm]{Definition}
\newtheorem{rem}[thm]{Remark}
\newtheorem{example}[thm]{Example}

\def\qqq{\,,\quad~\forall}

\def\GL{{\rm GL}}

\def\Hom{{\rm Hom}}

\def\Ker{{\rm Ker}}

\def\Res{{\rm Res}}

\def\Spec{{\rm Spec\,}}
\def\Sp{{\rm Spec}\,}

\def\A{{\mathbb A}}
\def\C{{\mathbb C}}
\def\F{{\mathbb F}}
\def\K{{\mathbb K}}

\def\N{{\mathbb N}}
\def\P{{\mathbb P}}
\def\Q{{\mathbb Q}}
\def\R{{\mathbb R}}
\def\Z{{\mathbb Z}}

\def\cC{{\mathcal C}}
\def\cD{{\mathcal D}}

\def\cG{{\mathcal G}}

\def\cF{{\mathcal F}}

\def\cO{{\mathcal O}}

\def\cS{{\mathcal S}}

\def\cT{{\mathcal T}}
\def\cU{{\mathcal U}}

\newcommand{\ie}{{\it i.e.\/}\ }

\newcommand{\cf}{{\it cf.\/}\ }
\newcommand{\opcit}{{\it op.cit.\/}\ }

\def\Hom {{\mbox{Hom}}}

\def\Gm{{\mathbb G}_m}

\def\ffp{\mathfrak{p}}

\def\strong{{\cS_\infty}}

\def\Mo{\mathfrak{Mo}}
\def\An{\mathfrak{ Ring}}
\def\Ab{\mathfrak{ Ab}}
\def\Mr{\mathfrak{ MR}}

\def\Se{\mathfrak{ Sets}}
\def\ssp{{\mathfrak{ spec}\,}}

\def\rep{\vartheta}
\def\spad{{\P^1_{\F_1}}}
\def\spadu{{\underline\P^1_{\F_1}}}
\def\Omm{\Omega}

\newcommand{\nil}[1]{}

\title
{Schemes over $\F_1$ and zeta functions}

\author[Connes]{Alain Connes}
\author[Consani]{Caterina Consani}
\address{A.~Connes: Coll\`ege de France \\
3, rue d'Ulm \\ Paris, F-75005 France
\\ I.H.E.S. and Vanderbilt
University} \email{alain\@@connes.org}
\address{C.~Consani: Mathematics Department \\ Johns Hopkins
University \\ Baltimore, MD 21218 USA} \email{kc\@@math.jhu.edu}
\thanks{The authors are partially supported by the NSF grant
DMS-FRG-0652164.}

\begin{document}
\maketitle

\begin{abstract} We  determine the {\em real} counting function $N(q)$ ($q\in [1,\infty)$) for the hypothetical ``curve'' $C=\overline {\Sp\,\Z}$ over $\F_1$, whose corresponding zeta function is the complete Riemann zeta function. Then, we develop a theory of functorial $\F_1$-schemes  which reconciles the previous attempts by C.~Soul\'e and A.~Deitmar. Our construction fits with the geometry of monoids of K.~Kato, is no longer limited to toric varieties  and it covers the case of schemes associated to Chevalley groups. Finally we show, using the monoid of ad\`ele classes over an arbitrary global field,
how to apply our functorial  theory of $\Mo$-schemes  to interpret conceptually the spectral realization of zeros of $L$-functions.
\end{abstract}

 \tableofcontents

\section{Introduction}

In this paper we develop three correlated aspects pertaining to the broad theory of the ``field of characteristic one'': $\F_1$. The appearance, in printed literature,  of some explicit remarks related to  this (hypothetical) degenerate algebraic structure is due to J.~Tits, who proposed its existence  to explain the limit case of the algebraic structure underlying the geometry of a Chevalley group over a finite field $\F_q$, as $q$ tends to $1$ (\cf~\cite{Tits}, \S~13 and \cite{ak}). A suggestive comment pointing out to a finite geometry inherent to the limit case $q=1$  is also contained in an earlier paper by R.~Steinberg (\cf~\cite{Steinberg}, p.~279), in relation to a geometric study of the representation theory of the general linear group over a finite field.

 In more recent years, the classical point of view that adjoining roots of unity is analogous to producing extensions of a base field, has  also been applied in the process of developing a suitable arithmetic theory  over $\F_1$. This idea lead to the introduction of the notion of algebraic field extensions $\F_{1^n}$ of $\F_1$ which are not defined {\it per se}, but are described by the following equation (\cf~\cite{Soule}, \S~2.4 and \cite{Kapranov})
  \[\label{basic}
\F_{1^n}\otimes_{\F_1} \Z:=\Z[T]/(T^n-1)\,, \qquad n\in\N.
\]
  The need for a field of characteristic one has also emerged in Arakelov's geometry, especially in the context of an absolute motivic interpretation of the zeros of zeta and L-functions (\cf~\cite{Manin}, \S~1.5). In \cite{Soule} (\S~6), C. Soul\'e introduced the zeta function of a variety $X$ over $\F_1$ by considering the {\em polynomial integer} counting function  of the associated  functor $\underline{X}$.\vspace{.05in}

 In this paper we take-up the following central question formulated in \cite{Manin} which originally motivated the development of the study of the arithmetic over $\F_1$.\vspace{.05in}

 {\bf Question}:~Can one find a ``curve'' $C=\overline {\Sp\,\Z}$ over $\F_1$ (defined in a suitable sense) whose zeta function $\zeta_C(s)$ is the complete Riemann zeta function  $\zeta_\Q(s)=\pi^{-s/2}\Gamma(s/2)\zeta(s)$? \vspace{.05in}

 After transforming the limit definition (for $q\to 1$) of the zeta function given in \cite{Soule} into an integral formula which is more suitable in the cases of general types of counting functions and distributions, we show how to determine the {\em real} counting function $N_C(q)=N(q)$, $q\in [1,\infty)$ associated to the hypothetical curve $C$ over $\F_1$.

 A convincing solution to this problem is a fundamental preliminary test for any  arithmetic theory over $\F_1$.
The difficulty inherent to the above question can be easily understood by considering the following facts. First of all, notice that the value $N(1)$ is conjectured to take the meaning of the Euler characteristic of the curve $C=\overline {\Sp\,\Z}$. Since one expects $C$ to be of infinite genus (\cf~\cite{Manin}), $N(1)$ is supposed to take the value  $-\infty$, thus precluding any easy use of the limit definition of the zeta and a naive approach to the definition of $C$, by generalizing the constructions of \cite{Soule}. On the other hand, the counting function $N(q)$ is also supposed to be positive for $q$ real, $q>1$, since it should detect the cardinality of the set of points of $C$ defined over various ``field extensions'' of $\F_1$. This requirement creates an apparent contradiction with the earlier condition $N(1)=-\infty$.

The precise statement of our result (\cf~Theorem \ref{fine1} and Remark \ref{Eulerchar}) is as follows:\newpage

$(1)$~The counting function $N(q)$ satisfying the above requirements exists and is given by the formula
  \begin{equation}\label{inghamn}
N(q)=q-\frac{d}{dq}\left(\sum_{\rho\in Z}{\rm order}(\rho)\frac{q^{\rho+1}}{\rho+1}\right)+1
\end{equation}
 where $Z$ is the set of non-trivial zeros of the Riemann zeta function and the derivative is taken in the sense of {\em distributions}.\vspace{.05in}

$(2)$~The function $N(q)$ is positive (as a {\em distribution}) for $q>1$.\vspace{.05in}

$(3)$~The value $N(1)$ is equal to $-\infty$ and reflects precisely the distribution of the   zeros of zeta in $E\log E$\footnote{$E=\frac{1}{\epsilon}$, $\epsilon>0$ appears when taking the derivative $\lim_{\epsilon\to 0} \frac{J(1+\epsilon)-J(1)}{\epsilon}$ of the primitive $J$ of $N(q)$}.\vspace{.05in}

This result supplies a strong indication on the coherence of the quest for an arithmetic theory over $\F_1$. Notice also that \eqref{inghamn} is  entirely similar to the classical formula for the counting function of the number of points of a proper, algebraic curve $X$ over $\F_p$ in the form
 $$
\#X(\F_q)=N(q)=q-\sum_{\alpha\in Z} \alpha^\ell+1,\qquad  \qqq q=p^\ell
$$
where the $\alpha$'s are the eigenvalues of the Frobenius operator acting on the cohomology of the curve.

The equation \eqref{inghamn} is a typical application of the Riemann-Weil explicit formulae. These formulae become natural when lifted to the id\`ele class group. This fact  supports the  expectation that, even if a definition of the hypothetical curve $C$ is at this time still out of reach, its counterpart, through the application of the class-field theory isomorphism, can be realized by a space of adelic nature and this is in agreement with  earlier constructions: \cf~\cite{CCM}, \cite{CCM2} and \cite{ccm}.\vspace{.05in}

A second topic that we develop in this paper is centered on the definition of a suitable geometric theory of algebraic schemes over $\F_1$. The viewpoint that we introduce in this article is an attempt at unifying the theories developed on the one side by Soul\'e  in \cite{Soule} and in our paper \cite{ak} and on the other side by A.~Deitmar in \cite{deit}, \cite{deit1}  (following  N.~Kurokawa,  H.~Ochiai and M.~Wakayama \cite{KOW}), by K.~Kato in \cite{Kato} (with the geometry of logarithmic structures) and by B.~T\"oen and  M.~Vaqui\'e in \cite{TV}.

In \cite{ak}, we introduced a refinement of the original notion (\cf~\cite{Soule}) of an affine variety over $\F_1$  and following this path we proved that Chevalley group schemes are examples of affine varieties defined over the field extension $\F_{1^2}$. While in the process of assembling this construction, we  realized that the functors (from finite abelian groups to graded sets)  describing these affine schemes fulfill stronger properties than the ones required in the original definition of Soul\'e. In this paper we develop this approach and show that the functors underlying the structure of the most common examples of schemes (of finite type) over $\F_1$  extend from (finite) abelian groups to a larger category obtained by gluing together the category $\Mo$ of commutative monoids (used in \cite{KOW}, \cite{Kato}, \cite{deit}, \cite{TV}) with the category  $\An$ of commutative rings. This  process uses a natural pair of adjoint functors relating $\Mo$ to $\An$ and follows an idea we learnt from P.~Cartier. The resulting category $\Mr$ (\cf~\S\ref{catmr} for details) defines an ideal framework in which the above two approaches are combined together to determine a very natural notion of variety (and of scheme) $\mathcal X$ over $\F_1$. In particular, the conditions imposed in the original definition of a variety over $\F_1$ in \cite{Soule} are now applied to a covariant functor $\mathcal X:\Mr \to \Se$ to the category of sets. Such a functor determines a scheme (of finite type) over $\F_1$ if it also fulfills the following three properties (\cf~Definition~\ref{defnfonesch}):\vspace{.05in}

-~The restriction $X_\Z$ of $\mathcal X$ to $\An$ is a scheme in the sense of \cite{demgab}.\vspace{.05in}

-~ The restriction $\underline X$ of $\mathcal X$ to $\Mo$ is locally representable.\vspace{.05in}

-~ The natural transformation connecting $\underline X$ to $X_\Z$, when applied to a field, yields a bijection (of sets).\vspace{.05in}

The category $\Ab$ of abelian groups embeds as a full subcategory in  $\Mo$. This fact allows one, in particular, to restrict a covariant functor from  $\Mo$ to sets to the subcategory  (isomorphic to) $\Ab$. In \S\ref{resabel} we prove that if the $\Mo$-functor is locally representable, then the restriction to $\Ab$ yields a functor to {\em graded} sets. This result shows that the grading structure that we assumed in \cite{ak} is now derived as a byproduct of this new refined approach.

In particular, we deduce that Chevalley groups are $\F_{1^2}$-schemes in our new sense; the group law exists on the set of points of {\em lowest degree} and is given by Tits' functorial construction of the normalizer of a maximal split torus.

As an arithmetic application of our new theory of $\F_1$-schemes we compute the zeta function of a Noetherian $\F_1$-scheme $\mathcal X$. Theorem~\ref{dthmfonesch} extends Theorem~1 of \cite{deit2} beyond the toric case and states, under a local torsion free hypothesis on the scheme, the following results:\vspace{.05in}

(a) There exists a polynomial $N(x+1)$ with positive integral coefficients such that
$$
  \# \,\underline X(\F_{1^n})=N(n+1)\qquad  \forall~n\in \N.
  $$
(b)
For each finite field $\F_q$, the cardinality of the set of points of the $\Z$-scheme
  $X_\Z$ which are rational over $\F_q$ is equal to $N(q)$.\vspace{.05in}

(c) The zeta function of $\mathcal X$ in the sense of \cite{Soule} is given by
  $$
    \zeta_\mathcal X(s)=\prod_{x\in X} \frac{1}{\left(1-\frac 1 s\right)^{\otimes^{n(x)}}}
$$
where the $\otimes$-product is the Kurokawa tensor product and  $n(x)$ denotes the local dimension at the point $x\in X$ of the geometric realization  of $\underline X$ (\cf~Definition~\ref{localdim}).\vspace{.05in}

The geometric theory of schemes over $\F_1$ that we have developed in \S\S~\ref{sectMo} and \ref{catmr} also reveals the importance to replace, when necessary, an abelian group $H$ by a naturally associated commutative monoid $M$ (with a zero element) so that $H=M^\times$ is interpreted as the group of invertible elements in the monoid. This idea applies in particular to the id\`ele class group $C_\K$ of a global field $\K$, since by construction the group $C_\K$ is the group of invertible elements in the multiplicative monoid of the ad\`ele classes
\begin{equation}\label{monoid}
    M=\A_\K/\K^\times\,, \ \K^\times=\GL_1(\K).
\end{equation}
 This application of the theory of $\Mo$-schemes to the study of geometric objects more pertinent to the realm of noncommutative geometry determines the third aspect of the theory of $\F_1$ that we have developed in this paper.  In our previous work, the ad\`ele class space has been considered mostly as a noncommutative space and its algebraic structure as a monoid did not play any role. One of the goals  of the present paper is to promote this additional structure by pointing out how and where it provides a precious guide.

 In \S\ref{projadel}, we consider  the particular case of the $\Mo$-scheme $\spadu$ describing a projective line over $\F_1$. It turns out that this scheme provides a perfect geometric framework to understand simultaneously and at a conceptual level, the spectral realization of zeros of $L$-functions, the functional equation and the explicit formulae. All these statements are deduced by simply computing the cohomology of a natural sheaf $\Omm$ of functions on the set $\spadu(M)$.  The  sheaf  $\Omm$ is a sheaf of complex vector spaces over the geometric realization $\spad$ of the $\Mo$-scheme $\spadu$. To define it we use a specific property of an $\Mo$-scheme, namely the existence, for each monoid $M$, of a natural projection $\pi_M:\underline X(M)\to X$, connecting the $\Mo$-scheme $\underline X$  (understood as  a functor from the category $\Mo$ of monoids to sets) to its associated geometric space $X$, \ie its geometric realization. For the $\Mo$-scheme $\spadu$ the geometric realization $\spad$ is a very simple space (\cite{deit}) which consists of three points
 $$
 \spad=\{0,u,\infty\}\, , \ \ \overline{\{0\}}=\{0\}\, , \ \overline{\{u\}}=\spad\, , \ \ \overline{\{\infty\}}=\{\infty\}.
$$
A striking fact is that in spite of the apparent simplicity of this space, the computation of $H^0(\spad,\Omm)$ already yields the graph of the Fourier transform: \cf~Lemma~\ref{lemsurject}. While the Fourier transform at the level of the ad\`eles depends upon the choice of a basic character, this dependence disappears at the level of the quotient space $M$ of ad\`ele classes. Also we explicitly remark that while the singularity of the operation $x\mapsto x^{-1}$ on the space of ad\`eles  prevents one to obtain any interesting global function on the projective space of the ad\`eles, this difficulty disappears at the level of the quotient space $M$ of ad\`ele classes (in view of the above result on $H^0(\spad,\Omm)$).

Theorem \ref{lemreasspec} states that the first cohomology group $H^1(\spad,\Omm)$ of the sheaf $\Omm$ over $\spad$ of complex valued functions on the projective space $\spadu(M)$ provides the space of the spectral realization of the zeros of $L$-functions. The symmetry associated to the functional equation derives as a simple consequence of the inversion $x\mapsto x^{-1}$ holding on $\spad$.

Finally, we want to stress the point that the most interesting aspect of this final  result does not rely on its technical part, since for instance the afore mentioned spectral realization is identical to that obtained in several earlier works  \cf~\cite{Meyer}, \cite{CCM}, \cite{CMbook} and initiated in \cite{Co-zeta}. The novelty of our statement is  that of proposing a new conceptual explanation for some fundamental constructions of noncommutative arithmetic geometry, in a way that the Fourier transform, the Poisson formula and the cokernel of the restriction map to the id\`eles all appear in an effortless and natural manner
on the projective line $\spadu(M)$.

\section{Zeta functions over $\F_1$ and $C=\overline {\Sp\,\Z}$}

In \cite{Soule} (\cf~\S 6) C. Soul\'e introduced the zeta function of a variety $X$ over $\F_1$ using the {\em polynomial} counting function $N(x)\in\Z[x]$ of the associated  functor $\underline{X}$. After correcting a sign misprint (which is  faithfully reproduced in \cite{deit2}), the precise definition of the zeta function is as follows
\begin{equation}\label{zetadefn}
\zeta_X(s):=\lim_{q\to 1}Z(X,q^{-s}) (q-1)^{N(1)},\qquad s\in\R
\end{equation}
where $Z(X,q^{-s})$ denotes the evaluation at $T=q^{-s}$  of the Hasse-Weil exponential series
\begin{equation}\label{zetadefn1}
Z(X,T) := \exp\left(\sum_{r\ge 1}N(q^r)\frac{T^r}{r}\right).
\end{equation}
Notice, incidentally,  that so defined $\zeta_X(s)$ fulfills the properties of an absolute motivic zeta function as predicted by Y. Manin in  \cite{Manin} (\cf~\S 1.5).

In this section, after transforming the limit \eqref{zetadefn} into an integral formula which is more suitable when dealing with general counting functions and distributions, we shall determine a precise formula for the counting function associated to the hypothetical curve $C=\overline{\Sp(\Z)}$.\vspace{.05in}

\subsection{An integral formula for  $\displaystyle{\partial_s\zeta_N(s)/\zeta_N(s)}$}\hfill\label{integral}\vspace{.05in}

Let $N(q)$ be a {\em real} continuous function on $[1,\infty)$ satisfying a polynomial bound $|N(q)|\leq C q^k$, for some finite positive integer $k$ and a fixed positive   constant $C$. Then, the corresponding generating function takes the following form
$$
Z(q,T)={\rm exp}\left( \sum_{r\geq 1}N(q^r)\frac{T^r}{r} \right)
$$
and one knows that the power series $Z(q,q^{-s})$ converges for $\Re(s)>k$. The zeta function over $\F_1$ associated to $N(q)$ is
$$
\zeta_N(s):=\lim_{q\to 1}Z(q,q^{-s})(q-1)^\chi \,, \ \ \chi=N(1)
$$
and this definition requires some care to assure its convergence. To eliminate the ambiguity in the extraction of the finite part, one works  with the logarithmic derivative
\begin{equation}\label{normael}
    \frac{\partial_s\zeta_N(s)}{\zeta_N(s)}=-\lim_{q\to 1} F(q,s)
\end{equation}
where
\begin{equation}\label{fqsdefn}
F(q,s)=-\partial_s \sum_{r\ge 1}N(q^r)\frac{q^{-rs}}{r}\,.
\end{equation}

\begin{lem} \label{compute1}With the above notations and for $\Re(s)>k$, one has
\begin{equation}\label{lim}
    \lim_{q\to 1} F(q,s) = \int_1^\infty N(u)u^{-s}d^*u\,,\ \ d^*u=du/u
\end{equation}
and
\begin{equation}\label{logzetabis}
    \frac{\partial_s\zeta_N(s)}{\zeta_N(s)}=-\int_1^\infty  N(u)\, u^{-s}d^*u\,.
\end{equation}
\end{lem}

\proof The proof follows immediately by noticing that
$$
F(q,s)=\sum_{r\ge 1}N(q^r) \,q^{-rs}\log q
$$
is the Riemann sum for the integral $\int_1^\infty N(u)u^{-s}d^*u$.
\endproof

Let us first assume  that $N(1)=0$. We then get the following expression by integrating in $s$, with $c$ a constant of integration,
\begin{equation}\label{logzeta}
    \log(\zeta_N(s))=\int_1^\infty \frac{N(u)}{\log u}u^{-s}d^*u + c.
\end{equation}
In the general case (\ie when $N(1)\neq 0$) one has to choose a principal value in the expression \eqref{logzeta} near $u=1$, since the term $\displaystyle{\frac{N(u)}{\log u}}$ is singular. The normalization used in  \cite{Soule}, corresponds to the principal value
\begin{equation}\label{logzetaprinc}
    \log(\zeta_N(s))=\lim_{\epsilon\to 0}\left(\int_{1+\epsilon}^\infty \frac{N(u)}{\log u}u^{-s}d^*u+N(1)\log \epsilon\right).
\end{equation}
Notice that this choice does not alter \eqref{logzetabis}. This fact is quite important since we shall use \eqref{logzetabis} to investigate the analytic nature of  $\zeta_N(s)$.\vspace{.05in}

\subsection{The counting function of $C=\overline {\Sp\,\Z}$}\hfill \label{distrsubsect}\vspace{.05in}

It is natural to wonder on the existence of a ``curve'' $C=\overline {\Sp\,\Z}$ suitably defined over $\F_1$, whose zeta function $\zeta_C(s)$ is the complete Riemann zeta function  $\zeta_\Q(s)=\pi^{-s/2}\Gamma(s/2)\zeta(s)$ (\cf also \cite{Manin}). In this subsection we shall show that
the integral equation \eqref{logzetabis} produces a precise formula for the counting function $N_C(q)=N(q)$ associated to $C$. In fact, \eqref{logzetabis} shows in this case that
\begin{equation}\label{special}
   \frac{\partial_s\zeta_\Q(s)}{\zeta_\Q(s)}=-\int_1^\infty  N(u)\, u^{-s}d^*u\,.
\end{equation}
This integral formula appears  in the Riemann-Weil explicit formulae and when $\Re(s)>1$, one derives that
\begin{equation}\label{special1}
   -\frac{\partial_s\zeta_\Q(s)}{\zeta_\Q(s)}=\sum_{n=1}^\infty\Lambda(n)n^{-s}+\int_1^\infty  \kappa(u)\, u^{-s}d^*u\,,
\end{equation}
where $\Lambda(n)$ is the von-Mangoldt function\footnote{with value $\log p$ for powers $p^\ell$ of  primes and zero otherwise} and $\kappa(u)$ is the distribution
$$
\kappa(u)=\frac{u^2}{u^2-1}\qquad  \forall~u>1
$$
which is defined using a principal value to eliminate the divergence at $u=1$. More precisely, the distribution $\kappa(u)$ is defined, for any test function $f$, by
$$
\int_1^\infty\kappa(u)f(u)d^*u=\int_1^\infty\frac{u^2f(u)-f(1)}{u^2-1}d^*u+cf(1)\,, \qquad c=\frac12(\log\pi+\gamma)
$$
where $\gamma=-\Gamma'(1)$ is the Euler constant.
Hence, we derive the consequence that the counting function $N(q)$ of the hypothetical curve $C$ over $\F_1$, is the {\em distribution} given by the sum of $\kappa(q)$ with the discrete term equal to the derivative $\frac{d}{dq}\varphi(q)$, taken in the sense of distributions, of the function\footnote{the value at the points of discontinuity does not affect the distribution}
\begin{equation}\label{varphifunc}
    \varphi(u)=\sum_{n<u}n\,\Lambda(n).
\end{equation}
Indeed, since  $d^*u=\frac{du}{u}$, one can write \eqref{special1} as
\begin{equation}\label{special2}
   -\frac{\partial_s\zeta_\Q(s)}{\zeta_\Q(s)}=\int_1^\infty \left(\frac{d}{du}\varphi(u)+ \kappa(u)\right)\, u^{-s}d^*u.
\end{equation}
If one compares \eqref{special2} and \eqref{special}, one derives the following formula for $N(u)$
\begin{equation}\label{Nu}
    N(u)=\frac{d}{du}\varphi(u)+ \kappa(u).
\end{equation}
The above expression encloses in a very subtle and intrinsic form a fundamental information on the description of the counting function as geometric ``trace type'' formula. To substantiate this statement, we recall the well-known equation (\cf~\cite{Ingham}, Chapter IV, Theorems 28 and 29, and use $\varphi(u)=u\psi_0(u)-\psi_1(u)$, where $\psi_0(u)$ is the Chebyshev function $\psi_0(u)=\sum_{n<u}\Lambda(n)$ and $\psi_1(u)=\int_0^u \psi_0(x)dx$ is its primitive) valid for $u>1$ (and not a prime power)
\begin{equation}\label{ingham1}
\varphi(u)=\frac{u^2}{2}-\sum_{\rho\in Z}{\rm order}(\rho)\frac{u^{\rho+1}}{\rho+1}+a(u)
\end{equation}
where
$$
a(u)={\rm ArcTanh}(\frac 1u)- \frac{\zeta'(-1)}{\zeta(-1)}
$$
  and $Z$ denotes the set of non-trivial zeros of the Riemann zeta function. Notice that  the sum over $Z$ in \eqref{ingham1} has to be taken in a symmetric manner to ensure convergence. When one differentiates it in a formal way, the term in $a(u)$ gives
$$
\frac{d}{du}a(u)=\frac{1}{1-u^2}.
$$
Hence, at the formal level \ie disregarding the principal value, one obtains
$$
\frac{d}{du}a(u)+\kappa(u)=1.
$$
Thus, when one differentiates (at the formal level) \eqref{ingham1} one gets
\begin{equation}\label{formallevel}
    N(u)=\frac{d}{du}\varphi(u)+ \kappa(u)\sim u-\sum_{\rho\in Z}{\rm order}(\rho)\,u^{\rho}+1.
\end{equation}
This formula for the counting function is now entirely similar to that describing the counting function of the number of points of a curve $C$ over a finite prime field $\F_p$ in the form of
 $$
\#C(\F_q)=N(q)=q-\sum_{\alpha\in Z} \alpha^\ell+1, \qquad\forall~q=p^\ell
$$
where the $\alpha$'s are the eigenvalues of the Frobenius.

Notice that in the above formal computations we have neglected to consider the principal value for the distribution $\kappa(u)$. By taking this into account we obtain the more precise result
\begin{thm}\label{fine1}
The distribution $N(u)$ satisfying the equation
$$
-\frac{\partial_s\zeta_\Q(s)}{\zeta_\Q(s)}=\int_1^\infty  N(u)\, u^{-s}d^*u\,,
$$
is positive on $(1,\infty)$ and  is given on $[1,\infty)$ by
\begin{equation}\label{fin2}
    N(u)=u-\frac{d}{du}\left(\sum_{\rho\in Z}{\rm order}(\rho)\frac{u^{\rho+1}}{\rho+1}\right)+1
\end{equation}
where the derivative is taken in the sense of distributions, and the value at $u=1$ of the  term
 $\displaystyle{\omega(u)=\sum_{\rho\in Z}{\rm order}(\rho)\frac{u^{\rho+1}}{\rho+1}}$ is given  by
\begin{equation}\label{definethet}
\omega(1)=\frac 12+ \frac \gamma 2+\frac{\log4\pi}{2}-\frac{\zeta'(-1)}{\zeta(-1)}.
\end{equation}
\end{thm}

\proof The positivity of the distribution $N(u)$ on $(1,\infty)$ follows from \eqref{Nu}. For $u>1$ we define
\begin{equation}\label{omega}
\omega(u)=\sum_{\rho\in Z}{\rm order}(\rho)\frac{u^{\rho+1}}{\rho+1}.
\end{equation}
By \eqref{ingham1} one has (for $u>1$)
\begin{equation}\label{omegabis}
\omega(u)=-\varphi(u)+\frac{u^2}{2}+a(u).
\end{equation}
In a neighborhood of $1$ one has $\varphi(u)=0$ and $a(u)\sim -\frac12\log(u-1)$ when $u\to 1+$.
Thus $\omega(u)$ diverges when $u\to 1$ although it is locally integrable and defines a distribution. Since $[1,\infty)$ has a boundary, the derivative of the distribution depends on its boundary value and is defined, for $f$ smooth and of fast enough decay at $\infty$, as
\begin{equation}\label{boundary}
\langle \frac{d}{du}\omega(u),f(u)\rangle=-\int_1^\infty\omega(u)\frac{d}{du}f(u)du-\omega(1)f(1).
\end{equation}
 We apply this to the function $f(u)=u^{-s-1}$, for $\Re(s)>1$.
One has $-\frac{d}{du}f(u)=(s+1) u^{-s-2}$ and one obtains
$$
\langle \frac{d}{du}\omega(u),f(u)\rangle=(s+1)\int_1^\infty
\left(-\varphi(u)+\frac{u^2}{2}+a(u)\right)u^{-s-2}du-\omega(1).
$$
By applying some results from \cite{Ingham} (\cf~Chapter~I, (17): use $\varphi(u)=u\psi_0(u)-\psi_1(u)$, $\psi'(u)=\psi_0(u)$), one deduces
\begin{equation}\label{phiphi}
-\frac{\partial_s\zeta(s)}{\zeta(s)}=(s+1)\int_1^\infty \varphi(u)u^{-s-2}du
\end{equation}
 and by using
$$
\int_1^\infty(u+1)f(u)du=\frac 1s+\frac{1}{s-1},\qquad -(s+1)\int_1^\infty \frac{u^2}{2}u^{-s-2}du=
-\frac 12-\frac{1}{s-1}
$$
one concludes that
$$
\langle (u-\frac{d}{du}\omega(u)+1),f(u)\rangle=\frac 1s-\frac 12-\frac{\partial_s\zeta(s)}{\zeta(s)}+\omega(1)-(s+1)\int_1^\infty
a(u)u^{-s-2}du.
$$
Finally, we claim that the following equation holds
$$
\frac 1s-(s+1)\int_1^\infty
a(u)u^{-s-2}du=-\frac{\partial_s\Gamma(s/2)}{\Gamma(s/2)}+ \frac{\zeta'(-1)}{\zeta(-1)}-\log 2-\frac \gamma 2.
$$
\begin{figure}
\begin{center}
\includegraphics[scale=1.4]{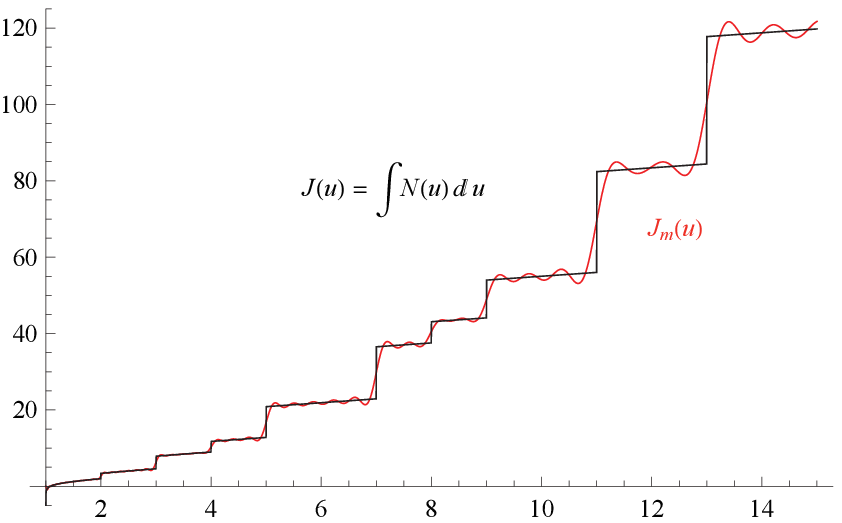}
\end{center}
\caption{{\it Primitive $J(u)$ of  $N(u)$ and approximation using the symmetric set $Z_m$ of first $2m$ zeros, by $$J_m(u)=\frac{u^2}{2}-\sum_{Z_m}{\rm order}(\rho)\frac{u^{\rho+1}}{\rho+1}+u$$ Note that $J(u)\to -\infty$ when $u\to 1+$}. \label{plotf} }
\end{figure}
Indeed, using a process of integration by parts one has
$$
-(s+1)\int_{1+\epsilon}^\infty
\left({\rm ArcTanh}(\frac 1u)-u\right)u^{-s-2}du=\int_{1+\epsilon}^\infty\frac{u^{-s+1}}{u^2-1}du+b(\epsilon)$$ with
$$b(\epsilon)=-\left({\rm ArcTanh}(\frac{ 1}{1+\epsilon})-(1+\epsilon)\right)(1+\epsilon)^{-s-1}=
-\int_{1+\epsilon}^\infty\frac{u^{-1}}{u^2-1}du+c+O(\epsilon \log(1/\epsilon))
$$
and $c=1-\log 2$. Moreover one also knows that
$$
\frac{\partial_s\Gamma(s/2)}{\Gamma(s/2)}=-\frac \gamma 2+ \int_1^\infty \frac{u^{-1}-u^{-s+1}}{u^2-1}du.
$$
Thus one gets
$$
\langle (u-\frac{d}{du}\omega(u)+1),f(u)\rangle=-\frac{\partial_s\zeta_\Q(s)}{\zeta_\Q(s)},
$$
provided that
\begin{equation}\label{omega1}
\omega(1)=\frac 12+ \frac \gamma 2+\frac{\log4\pi}{2}-\frac{\zeta'(-1)}{\zeta(-1)}.
\end{equation}
To check this latter equality one cannot use the explicit formula \eqref{ingham1} which is not valid at $u=1$, since the term ${\rm ArcTanh}(\frac 1u)$ is infinite, therefore displaying the discontinuity of the function $\omega(u)$ at $u=1$. To verify \eqref{omega1} we rather use the following formula (taken from \cite{Ingham}, \cf~III, (26))
$$
\frac{\zeta'(s)}{\zeta(s)}+\frac{1}{s-1}=\sum_Z\left(\frac{1}{s-\rho}+\frac{1}{\rho}\right)-\frac 12\frac{\Gamma'}{\Gamma}(\frac s2+1)+\log(2\pi)-1-\frac \gamma 2
$$
when $s\to 1$. We notice that the left hand side of the above formula tends to $\gamma$ while the right hand side, using the symmetry $\rho\to 1-\rho$ of the zeros (and a symmetric summation and the formula $
\frac{\Gamma'}{\Gamma}\left(\frac 32\right)=2-\gamma-2\log 2
$) tends to
$$
2 \sum_Z\frac{1}{\rho}-2+\log(4\pi).
$$
Thus one obtains
$$
\sum_Z\frac{1}{\rho}=\frac \gamma 2+1-\frac 12 \log(4\pi).
$$
One then concludes by using the equalities (\cf \cite{Ingham} IV, Theorem 28)
$$
\sum_Z\frac{1}{\rho(\rho+1)}=\frac 12-\log(4\pi)+\frac{\zeta'(-1)}{\zeta(-1)}
$$
and the formula (using a symmetric summation)
$$
\sum_Z\frac{1}{\rho+1}=\sum_Z\frac{1}{\rho}-\sum_Z\frac{1}{\rho(\rho+1)}.
$$
\endproof

\begin{rem}\label{Eulerchar}{\rm
In agreement with \cite{Soule}, the value $N(1)$ should be thought of as the Euler characteristic
of the hypothetical curve $C$ over $\F_1$. Since $C$ is expected to have infinite genus, one would deduce that $N(1)=-\infty$,  in apparent conflict with  the expected positivity of $N(q)$ for $q>1$. This apparent contradiction is resolved in the proof of Theorem~\ref{fine1}, since the distribution $N(q)$ is positive for $q>1$ but its value at $q=1$ is formally given by
$$
N(1)=2-\lim_{\epsilon\to 0}\frac{\omega(1+\epsilon)-\omega(1)}{\epsilon}\sim-\frac 12 E \log E,\qquad \ E=\frac 1\epsilon
$$
reflecting, when $\epsilon\to 0$, also the density of the zeros.
}\end{rem}

Equality \eqref{ingham1} is a typical application of the Riemann-Weil explicit formulae which become natural once they are lifted to the id\`ele class group. It seems therefore natural also to expect that the hypothetical curve $C=\overline{\Sp(\Z)}$ is of ad\`elic nature and that it also possesses an action of the id\`ele class group. This speculation is in agreement with the interpretation of the explicit formulae as a trace formula, by using the noncommutative geometric formalism of the ad\`ele class space (\cf \cite{Co-zeta},  \cite{CCM}, \cite{CCM2}, \cite{ccm}, \cite{Meyer}).

 \section{$\Mo$-schemes.}\label{sectMo}

 In this section we describe, following a functorial approach similar to that of \cite{demgab}, a generalization of the theory of $\Z$-functors and schemes obtained by enlarging the category of rings to that of commutative monoids. This functorial construction will be applied in \S\ref{funfunctor}, after gluing together the categories of monoids and rings, to derive a new notion of $\F_1$-schemes and associated zeta functions. Our construction has  evident connections with the theory of schemes over $\F_1$ developed by A.~Deitmar in \cite{deit}, \cite{deit1}, with the theory of logarithmic structures of K.~Kato in  \cite{Kato}, with the arithmetic theory over $\F_1$ described by N.~Kurokawa,  H.~Ochiai, M.~Wakayama in \cite{KOW}, and with the algebro-topological approach followed by B.~T\"oen and  M.~Vaqui\'e in \cite{TV}.\vspace{.05in}

 \subsection{Monoids: the category $\Mo$.}\label{monoids}\hfill\vspace{.05in}

Throughout the paper we denote by $\Se$, $\Ab$, $\An$ respectively the categories of sets, abelian groups and commutative rings with unit.

We let $\Mo$ be the category of commutative  monoids $M$ denoted multiplicatively, with a neutral element $1$ (\ie unit) {\em and} an absorbing element  $0$ ($0\cdot x =x\cdot 0= 0,~\forall x\in M$).

A homomorphism $\varphi: M\to N$ in $\Mo$ is  unital (\ie $\varphi(1) = 1$) and satisfying $\varphi(0) = 0$.

Given a commutative group $H$ in $\Ab$, we set
\[
\F_1[H] = H\cup\{0\}\qquad (0\cdot h = h\cdot 0 = 0,\quad\forall h\in H).
\]
 Following the analogy with the category of rings, one sees that in $\Mo$ a monoid of the form $\F_1[H]$ corresponds to a {\em field} $F$ ($F = F^\times\cup\{0\}$) in $\An$. The collection of monoids like $\F_1[H]$, for $H\in\text{Obj}(\Ab)$, forms a full subcategory of $\Mo$ isomorphic to the category of abelian groups: \cf~Proposition~\ref{emb}.

\begin{defn} \label{Mfunc} An $\Mo$-functor $\cF$ is a covariant functor from the category $\Mo$ to $\Se$.
 \end{defn}

 To a monoid $M$  in $\Mo$ one associates the covariant functor $\ssp M$ defined as follows
\begin{equation}\label{specm}
\ssp M: \Mo \to \Se\qquad N\mapsto \ssp M(N)=\Hom_{\Mo}(M,N).
\end{equation}
 Notice that by applying Yoneda's lemma, a morphism of functors (natural transformation)  $\varphi:\ssp M\to \cF$, with $\cF: \Mo \to \Se$ is completely determined by the element $\varphi(id_M)\in \cF(M)$, moreover any such element gives rise to a morphism $\ssp M\to \cF$. By applying this fact to the functor $\cF = \ssp N$, for $N\in\text{Obj}(\Mo)$, one obtains an {\em inclusion} of $\Mo$ as a full subcategory of the category of $\Mo$-functors.\vspace{.05in}

Morphisms in the category of $\Mo$-functors are natural transformations.\vspace{.05in}

 An ideal  $I$ of a monoid $M$ is a subset $I\subset M$ such that $0\in I$ and
 $x\in I \implies xy\in I\qqq y\in M $ (\cf \cite{Gilmer}).
 As for rings, an ideal $I\subset M$ defines an interesting subfunctor $\underline{D(I)}\subset\ssp M$:
 \begin{equation}\label{di}
 \underline{D(I)}: \Mo \to \Se,\qquad \underline{D(I)}(N) = \{\rho\in\ssp(M)(N)|\rho(I)N=N\}.
 \end{equation}
  We recall that an ideal $\ffp\subset M$ is said to be {\em prime} if $1\notin\ffp$ and its complement  $\ffp^c= M\setminus\ffp$ is a multiplicative subset of $M$ \ie
 $$x\notin \ffp,~y\notin \ffp \implies xy \notin \ffp.$$

For an ideal $I\subset M$, one denotes by $D(I)$  the set of prime ideals $\ffp\subset M$ which do not contain $I$. These subsets are the open sets for the natural topology on the set $X=\Spec(M)$ of prime ideals of $M$ (\cf~\cite{Kato}). The smallest ideal containing a collection of ideals $\{I_\alpha\}$ of a monoid $M$ is just the union $I = \cup_\alpha I_\alpha$ and the corresponding open subset $D(I)\subset \Sp (M)=\{\ffp\subset M|\ffp~\text{prime ideal}\}$  satisfies the property  $D(\cup_\alpha I_\alpha) = \cup_\alpha D(I_\alpha)$. It is a standard fact that the inverse image of a prime ideal by a morphism of monoids is a prime ideal. Moreover, it is also straightforward to verify that the complement of the set of invertible elements in a monoid $M$, $\ffp_M=(M^\times)^c$, is a prime ideal in $M$ which contains all other prime ideals of the monoid.\vspace{.05in}

 \subsection{Automatic locality.}\hfill\label{autloc}\vspace{.05in}

 An interesting property fulfilled by any $\Mo$-functor is that of being local. Locality is not automatically satisfied by $\Z$-functors, essentially it corresponds to state the exactness, on an open covering of an affine scheme $\Spec(R)=\cup_i D(f_i)$ ($f_i\in R$), of sequences such as \eqref{locf} below.
 On the other hand, we shall see that an $\Mo$-functor is local by construction.  We recall the following result (\cf~\cite{deit})

\begin{lem}\label{local}  Let $M$ be an object in $\Mo$ and let $\{W_\alpha\}_{\alpha\in A}$ ($A$ a set) be an open cover of  the topological space $X=\Spec(M)$. Then $W_\alpha = \Spec(M)$, for some index $\alpha\in A$.
\end{lem}
\proof The point $\ffp_M = (M^\times)^c\in\Spec(M)$ must be contained in at least one  $W_\alpha$, for some index $\alpha\in A$. One has $W_\alpha=D(I_{\alpha})$ for some ideal $I_\alpha\subset M$,  hence $\ffp_M \in D(I_{\alpha})$, for some $\alpha\in A$  and this means $I_{\alpha}\cap M^\times \neq\emptyset$, that is $I_{\alpha} = M$.
\endproof

 Let $M$ be an object of $\Mo$. For $S\subset M$ a {\em multiplicative} subset we recall that the {\em monoid} $S^{-1}M$ is the quotient of the set made by all  expressions $a/s=(a,s)\in A\times S$, by the following equivalence relation
\[
a/s\sim b/t\quad \Leftrightarrow\quad\exists~u\in S\quad ut a = usb.
\]
One checks that the product
$
a/s.b/t=ab/st
$
is well-defined on the quotient $S^{-1}M$. For $f\in M$ and $S=\{f^n; n\in \Z_{\ge 0}\}$ one denotes $S^{-1}M$ by $M_f$.

For any $\Mo$-functor $\cF$ and any monoid $M$  one defines a sequence of maps of sets
\begin{equation}\label{locf}
    \cF(M)\stackrel{u}{\longrightarrow}\prod_i\cF(M_{f_i})\;\xy
{\ar@{->}_{w} (0,-1)*{}; (6,-1)*{}};
{\ar@{->}^{v} (0,1)*{}; (6,1)*{}};
\endxy\prod_{ij}
    \cF(M_{f_if_j})
\end{equation}
 which is obtained by using the open covering of $\Spec(M)$ made by the open sets $D(f_iM)$ ($f_i\in M)$, the natural morphisms $M \to M_{f_i}$ and the functoriality of $\cF$.\vspace{.05in}

 The following lemma shows that any $\Mo$-functor is local
\begin{lem} \label{lemauto} For any $\Mo$-functor $\cF$ and any monoid $M$, the sequence \eqref{locf} is exact.
\end{lem}
\proof By Lemma~\ref{local}, there exists an index $i$ such that $f_i\in M^\times$. One may assume that $i = 1$. Then, the map $\rho_1:M\to M_{f_1}$ is invertible thus $u$ is injective. Let $(x_i)\in \prod_i\cF(M_{f_i})$ be a family, with $x_i\in \cF(M_{f_i})$ such that $(x_i)_{f_j} = (x_j)_{f_i}$, for all $i,j$. This gives in particular the equality between the image of $x_i\in\cF(M_{f_i})$ under the isomorphism  $\cF(\rho_{i1}): \cF(M_{f_i}) \to \cF(M_{f_if_1})$ and $\cF(\rho_{1i})(x_1) \in \cF(M_{f_1f_i})=\cF(M_{f_if_1})$. By writing $x_1 = \rho_1(x)$ one finds that $u(x)$ is equal to the family $(x_i)$.
\endproof

 \subsection{Open $\Mo$-subfunctors.}\hfill\label{openMo}\vspace{.05in}

 In analogy with the theory of $\Z$-schemes, we now introduce the notion of an open subfunctor

 \begin{defn} \label{opensubmo} A subfunctor $ \mathcal G\subset  \cF$ of an $\Mo$-functor $\cF$ is open if for any object $M$ of $\Mo$ and any morphism of $\Mo$-functors $\varphi: \ssp M \to  \cF$, there exists an ideal $I\subset M$ satisfying the following property\vspace{.05in}

For any object $N$  of $\Mo$ and for any $\rho\in\ssp M(N)=\Hom_{\Mo}(M,N)$:
\begin{equation}\label{condition1}
\varphi(\rho)\in  \mathcal G(N)\subset  \cF(N)~\Leftrightarrow~\rho(I)N = N.
\end{equation}
\end{defn}

To clarify the meaning of this definition we develop a few  examples.\vspace{.05in}

\begin{example} \label{mulsub} {\rm The functor
\[
 \mathcal G: \Mo\to \Se,\quad N \to  \mathcal G(N)=N^\times
 \]
 is an open subfunctor of the (identity) functor $\cD^1$
 \[
  \cD^1: \Mo \to \Se,\quad N \to  \cD^1(N)=N.
  \]
  In fact, let $M$ be a monoid, then by Yoneda's lemma a morphism of functors $\varphi: \ssp M \to  \cD^1$  is determined by an element $z\in \cD^1(M) = M$. For any monoid $N$ and $\rho\in\Hom(M,N)$, one has $\varphi(\rho)=\rho(z)\in  \cD^1(N)=N$, thus the condition $\varphi(\rho)\in  \mathcal G(N) = N^\times$ means that $\rho(z)\in N^\times$. One takes  for $I$ the ideal generated by $z$ in $M$: $I = zM$. Then it is straightforward to check that \eqref{condition1} is fulfilled.}\end{example}

  \begin{example} \label{opensubf} {\rm Let $I\subset M$ be an ideal of a monoid $M$ and consider the  subfunctor $\underline{D(I)}\subset\ssp(M)$ as defined in \eqref{di}.
  Then, $\underline{D(I)}$ is an open subfunctor of $\ssp M$.

  Indeed, for any object $A$ of $\Mo$ and $\varphi: \ssp A \to \ssp M$ one has $\varphi(id_A) = \eta\in\ssp(M)(A) = \text{Hom}_{\Mo}(M,A)$. One takes in $A$ the ideal $J = \eta(I)A$. This ideal fulfills the condition \eqref{condition1} for any object $N$ of $\Mo$ and $\rho\in\text{Hom}_{\Mo}(A,N)$.
In fact, one has $\varphi(\rho) = \rho\circ\eta\in\text{Hom}_{\Mo}(M,N)$ and $\varphi(\rho)\in \underline{D(I)}(N)$ means that $\rho(\eta(I))N = N$. This latter equality holds if and only if $\rho(J)N = N$.
  }\end{example}\vspace{.05in}

 \subsection{Open covering by $\Mo$-subfunctors.}\hfill\label{coverMo}\vspace{.05in}

 The next task is to introduce the notion of an {\em open cover} in the category of $\Mo$-functors.
  We shall use the fact (\cf~Proposition~\ref{emb}) that the category $\Ab$ of abelian groups embeds as a full subcategory of $\Mo$, by means of the functor $H\to \F_1[H]$.

\begin{defn} \label{defncover} Let $\cF$ be an $\Mo$-functor and let $\{{\cF_\alpha}\}_{\alpha\in S}$ be a family of open subfunctors of $\cF$. Then one says that $\{{\cF_\alpha}\}_{\alpha\in S}$ ($S=$ an index set) is an open cover of $\cF$ if
\begin{equation}\label{cover}
 \cF(\F_1[H]) = \bigcup_{\alpha\in S} {\cF_\alpha}(\F_1[H]),\quad\forall H\in{\rm Obj}(\Ab).
\end{equation}
\end{defn}
Since commutative groups (with $0$ added) replace fields in $\Mo$, the above definition is the natural transposition of the definition of open covers as in \cite{demgab} within the category of $\Mo$-functors. The following proposition gives a precise characterization of the open covers of an $\Mo$-functor

\begin{prop} \label{coverspec} Let $\cF$ be an $\Mo$-functor and let $\{{\cF_\alpha}\}_{\alpha\in S}$ be a family of open subfunctors of $\cF$. Then, the family $\{{\cF_\alpha}\}_{\alpha\in S}$ forms an open cover of $\cF$ if and only if  \[
\cF(M)=\bigcup_{\alpha\in S} \cF_\alpha (M),\quad\forall~M\in{\rm Obj}(\Mo).
\]
 \end{prop}

 \proof The condition is obviously sufficient. To show the converse, we assume \eqref{cover}. Let $M$ be a monoid and let $\xi\in \cF(M)$, one needs to show that $\xi\in \cF_\alpha (M)$ for some $\alpha\in S$.
 Let $\varphi$ be the morphism of functors from $\ssp M$ to $\cF$ such that $\varphi(id_M)=\xi$. Since each $\cF_\alpha$ is an open subfunctor of $\cF$, one can find ideals $I_\alpha\subset M$ such that for any object $N$  of $\Mo$ and for any $\rho\in\ssp M(N)=\Hom_{\Mo}(M,N)$ one has
 \begin{equation}\label{condcover}
 \varphi(\rho)\in  \cF_\alpha(N)\subset  \cF(N)~\Leftrightarrow~\rho(I_\alpha)N = N.
\end{equation}
One applies this to the morphism $\epsilon_M:M\to \F_1[M^\times]=\kappa$ given by
\begin{equation}\label{residuemorphism}
    M \stackrel{\epsilon}{\to}\F_1[M^\times]\,, \ \epsilon_M (y)=0\qqq y\notin M^\times\,, \ \epsilon_M (y)=y\qqq y\in M^\times.
\end{equation}
One has $\epsilon_M\in \ssp(M)(\kappa)$ and $\varphi(\epsilon_M)\in \cF(\kappa)=\bigcup_{\alpha\in S} {\cF_\alpha}(\kappa)$. Thus,
 $\exists \alpha$ such that  $\varphi(\epsilon_M)\in \cF_\alpha(\kappa)$. By \eqref{condcover} one has $\epsilon_M(I_\alpha)\kappa=\kappa$ and
$I_\alpha\cap M^\times\neq \emptyset$ hence $I_\alpha=M$. Applying then \eqref{condcover} to $\rho=id_M$ one obtains $\xi\in \cF_\alpha (M)$ as required. \endproof\vspace{.05in}

 \subsection{$\Mo$-schemes.}\hfill\label{Moscheme}\vspace{.05in}

 In view of the fact that any $\Mo$-functor is local, the definition of an
 $\Mo$-scheme simply involves the local representability

 \begin{defn}\label{defnmosch} An $\Mo$-scheme is an $\Mo$-functor which admits an open cover by representable subfunctors.
\end{defn}

We shall consider several elementary examples of $\Mo$-schemes\vspace{.05in}

\begin{example}\label{affineex}{\rm  The affine spaces $\cD^n$. For a fixed $n\in\N$, we consider the following $\Mo$-functor
\[
\cD^n: \Mo \to \Se,\quad \cD^n(M) = M^n
\]
This functor  is representable since it is described by
\[
\cD^n(M) =  \text{Hom}_{\Mo}(\F_1[T_1,\ldots ,T_n],M),
\]
where
 \begin{equation}\label{F_1T}
\F_1[T_1,\ldots,T_n] := \{0\}\cup\{T_1^{a_1}\cdots T_n^{a_n}|a_j\in\Z_{\ge 0}\}\,
\end{equation}
is the union of $\{0\}$ with the semi-group generated by the $T_j$.
}\end{example}

\begin{example}\label{projex}{\rm The projective line $\underline{\P^1}$.  We consider the $\Mo$-functor $\underline{\P^1}$ which associates to an object $M$ of $\Mo$  the set $\underline{\P^1}(M)$ of complemented submodules $E$ of rank one in $M^2$, where the rank is defined locally. By definition a complemented submodule is the range of an idempotent matrix $e\in M_2(M)$ (\ie $e^2=e$) with each line having at most\footnote{Note that we need the $0$-element to state this condition} one non-zero entry.   To a morphism $\rho:M\to N$ one associates the following map $\underline{\P^1}(\rho)$
$$
E\to N\otimes_M E\subset N^2
$$
which replaces  $e\in M_2(M)$ by $\rho(e)\in M_2(N)$. The condition of rank one means that for any prime ideal $\ffp\in \Sp M$ one has $\epsilon_\ffp(e)\notin\{0,1\}$ where $\epsilon_\ffp$ is the morphism from $M$ to $\F_1[M_\ffp^\times]$ (where  $M_\ffp=S^{-1}M$, with $S = \ffp^c$) given by
\begin{equation}\label{locloc}
   \epsilon_\ffp(y)=0\qqq y\in \ffp\,, \ \ \epsilon_\ffp(y)=y\qqq y\notin \ffp.
\end{equation}

 Now, we compare $\underline{\P^1}$ with the $\Mo$-functor
\begin{equation}\label{projlinemo}
\mathcal P(M) = M \cup_{M^\times}M
\end{equation}
where the gluing map is given by $x\to x^{-1}$. In other words, we define on the disjoint union $ M\cup M$ an equivalence relation given by (using the identification $M\times\{1,2\} = M\cup M$)
\[
(x,1)\sim (x^{-1},2)\quad\forall x\in M^\times.
\]
 We define a natural transformation $e$ from $\mathcal P$ to $\underline{\P^1}$ by observing that the matrices
$$
e_1(a)=\left(
  \begin{array}{cc}
    1 & 0 \\
    a & 0 \\
  \end{array}
\right)\,, \ \
e_2(b)=\left(
  \begin{array}{cc}
    0 & b \\
    0 & 1 \\
  \end{array}
\right)\,, \ \ a, b\in M
$$
are idempotent ($e^2=e$) and their ranges also fulfill the following property
$$
 {\rm Im}\,e_1(a)={\rm Im}\,e_2(b)\iff ab=1.
$$

\begin{lem}\label{projtwocomp} The natural transformation $e$ is an isomorphism \ie
$$
 \mathcal P(M)= M \cup_{M^\times}M \cong \underline{\P^1}(M).
$$
Moreover, the two copies of $M$ define an open cover of $\underline{\P^1}$ by representable sub-functors $\cD^1$.
\end{lem}
\proof We show that an idempotent matrix $e\in M_2(M)$ of rank one, with each line having at most one non-zero entry is of the form $e_j(a)$ for some $j\in\{1,2\}$.
First we claim  that one of the matrix elements of $$e=\left(
  \begin{array}{cc}
    a & b \\
    c & d \\
  \end{array}
\right)$$ must be invertible. Otherwise, by localizing $M$ at the prime ideal $\ffp_M=(M^\times)^c$ one would obtain the zero matrix  which contradicts the hypothesis of rank one.
Assume first that $a$ is invertible. Then $b=0$, and from the idempotency condition on $e$ one gets that $a^2=a$ and hence $a=1$. Now we show that $d=0$. Again from the condition $e^2=e$ one gets $d^2=d$. Then, if $d\neq 0$ there exists a prime ideal $\ffp\subset M$
such that $d\notin \ffp$. This because the intersection of all prime ideals is the set of nilpotent elements. More generally, one knows (\cite{Gilmer}) that given an ideal $I\subset M$, the intersection  of the prime ideals $\ffp\subset M$, with $\ffp\supset I$ coincides with the radical of $I$
\[
\bigcap_{\ffp\supset I}\ffp = \sqrt I := \{x\in M|\exists n\in\N, x^n\in I\}.
\]
Thus, by localizing $M$ at $\ffp$  one gets that $e$ is the unit matrix at $\ffp$ which contradicts the hypothesis of rank one.
Thus $d=0$ and $e=e_1(c)$. \vspace{.05in}

If $b$ is invertible then $a=0$, $bd=b$ so that $d=1$ and $c=0$ thus $e=e_2(b)$.
\vspace{.05in}

The two other cases are treated in a similar manner.
The functor $\mathcal P$ admits by construction two copies of the functor $\cD^1$ embedded in it as subfunctor. We need to show that these two subfunctors are open in $\underline{\P^1}$: we prove it for the first copy of $\cD^1$. Let $N$ be an object of $\Mo$. A morphism $\ssp (N) \to \underline{\P^1}$ (in the category of $\Mo$-functors) is determined by an element $z\in\underline{\P^1}(N)$. If $z$ belongs to the first copy of $N$, it follows that for any $\rho\in \text{Hom}_{\Mo}(N,M)$,  $\rho(z)$ is in the first copy of $\cD^1(M)$. In this case one can take $I = N$. Otherwise, $z$ belongs to the second copy of $N$ and in this case, likewise in the above Example \ref{mulsub}, one takes $I = zN$. The local representability follows since $\cD^1$ is representable.\endproof
}\end{example}

\begin{example}\label{opensubff}{\rm  Let $M$ be a monoid and let $I\subset M$ be an ideal. Consider the $\Mo$-functor $\underline{D(I)}$ of Example \ref{opensubf}. The next proposition states that this is an $\Mo$-scheme.
}
\end{example}

\begin{prop}\label{opensubfff} 1) Let $f\in M$ and $I = fM$. Then the subfunctor $\underline{D(I)}\subset\ssp M$ is represented by $M_f$.\vspace{.05in}

2) For any ideal $I\subset M$, the $\Mo$-functor
$\underline{D(I)}$ is an $\Mo$-scheme.
\end{prop}
The proof is straightforward.\vspace{.05in}

 \subsection{Geometric realization.}\hfill\label{geomereal}\vspace{.05in}

As in the case of $\Z$-schemes (and following similar set-theoretic precautions as the ones stated in the preliminary chapter of \cite{demgab}), it can be shown that any $\Mo$-scheme $\underline X$ can be represented in the form
 \begin{equation}\label{morep}
    \underline X(N)=\Hom(\Sp(N),X)\qquad N\in\text{Obj}(\Mo),
 \end{equation}
where $X$ is the associated geometric space, \ie the geometric realization of $\underline X$. In this framework, a geometric space is properly defined by\vspace{.05in}

$\bullet$~A topological space $X$\vspace{.05in}

$\bullet$~A sheaf $\cO_X$ of monoids on $X$.\vspace{.05in}

For details on the properties of the geometric spaces which are locally of the form $\Sp(M)$ we refer to  \cite{Kato}, \cite{deit} and \cite{deit1}. Notice that there is no need to require that the stalks of the structural sheaf of a geometric space are ``local'' since any monoid has already a local algebraic structure.
We recall only a few concepts from the basic terminology and we refer to  \cite{demgab}, \cite{Kato}, \cite{deit} and \cite{deit1} for details. A morphism  $\rho: M_1\to M_2$ of monoids is said to be local if
$
\rho^{-1}(M_2^*) = M_1^*.
$
A morphism $\varphi: X \to Y$ between two geometric spaces is given by a pair $(\varphi,\varphi^\sharp)$ of a continuous map $\varphi$ and a local morphism of sheaves of monoids
\[
\Gamma(V,\mathcal O_Y)\stackrel{\varphi^\sharp}{\to}\Gamma(\varphi^{-1}(V),\mathcal O_X)
\]
\ie  the map of stalks $\mathcal O_{\varphi(x)}\stackrel{\varphi^\sharp}{\to}\mathcal O_x$ is local.

The sheaf of monoids associated to the {\em prime spectrum} $\Sp(M)$  satisfies the following properties:\vspace{.05in}

$\bullet$~The stalk at $\ffp\in\Spec(M)$ is $\mathcal O_{\ffp}=S^{-1}M$,  with  $S = \ffp^c$.\vspace{.05in}

$\bullet$~For any $f\in M$, the  map $\varphi: M_f \to \Gamma(D(fM),\mathcal O)$ defined by
  \[
\varphi(x)(\ffp) = a/f^n\in \mathcal O_\ffp\quad \forall\ffp\in D(fM)\qqq x = a/f^n\in M_f
\]
is an isomorphism.\vspace{.05in}

$\bullet$~On an open set $U\subset\Spec(M)$, a section $s\in \Gamma(U,\mathcal O)$ is an element of $\prod_{\ffp\in U}\mathcal O_\ffp$ such that on any open set $D(f)\subset U$ its restriction agrees with an element in $M_f$.\vspace{.05in}

For any geometric space $(X,\mathcal O_X)$, one defines (as in \cite{demgab}) a canonical morphism $\psi_X:X\to \Sp(\cO(X))$.

\begin{defn} \label{geometricsch} A geometric space $(X,\mathcal O_X)$ is a {\em prime spectrum} if  the morphism $\psi_X$ is an isomorphism.
It is a {\em  geometric $\Mo$-scheme} if $X$ admits an open covering by prime spectra.
\end{defn}

The terminology is justified since the $\Mo$-functor $\underline X(M)=\Hom(\Sp M,X)$ associated to a geometric $\Mo$-scheme $X$ is an $\Mo$-scheme in the sense of Definition \ref{defnmosch}. We can now state
\begin{prop}\label{representability} Under the same set-theoretic conditions as in \cite{demgab}, any $\Mo$-scheme $\underline X$ can be represented in the form
 \begin{equation}\label{morep}
    \underline X(N)=\Hom(\Sp(N),X),
 \end{equation}
 for a geometric $\Mo$-scheme $X$ which is unique up to isomorphism.
\end{prop}\vspace{.05in}

The proof of this proposition follows the lines of that for $\Z$-schemes exposed in \cite{demgab}. The geometric realization $X$ ($|\underline X|$ in the notation of \cite{demgab}) is constructed canonically as an inductive limit of prime spectra: \cf Proposition 4.1 of \opcit Proposition \ref{coverspec} ensures that the natural map from $\Hom(\Sp(N),X)$ to $\underline X(N)$ is surjective.

\vspace{.05in}

 \subsection{Restriction to abelian groups.}\hfill \label{resabel}\vspace{.05in}

 In this section we describe the functor obtained by restricting $\Mo$-schemes to the category $\Ab$ of abelian groups. We first recall the definition of the natural functor-inclusion of $\Ab$ in $\Mo$.

 \begin{prop}\label{emb} The covariant functor
\[
\F_1[~\cdot~]: \Ab\to \Mo\qquad H\mapsto \F_1[H]
\]
embeds the category of abelian groups as a full subcategory of the category of commutative monoids.
\end{prop}
\proof We show that the group homomorphism
$$
\Hom_{\Ab}(H,K)\stackrel{}{\longrightarrow} \Hom_{\Mo}(\F_1[H],\F_1[K])\qquad \phi\to\F_1[\phi]
$$
is bijective. It is injective by restriction to $H\subset\F_1[H]$. Moreover, any unital monoid homomorphism in $\Hom_{\Mo}(\F_1[H],\F_1[K])$ preserves the absorbing elements and sends invertible elements to invertible elements since it is unital. Thus it arises from a group homomorphism.
\endproof

We shall identify $\Ab$ with this full subcategory of $\Mo$. Of course, any $\Mo$-functor $\underline X: \Mo\to \Se$ can be restricted to $\Ab$ and it gives rise to a functor taking values in $\Se$.

In fact, there is a pair of adjoint functors: $\Ab\to\Mo,~H\mapsto \F_1[H]$ and  $\Mo\to\Ab,~M\mapsto M^\times$  linked by the isomorphism
$$\Hom_{\Mo}(\F_1[H],M) \cong \Hom_{\Ab}(H,M^\times).$$

Moreover, for a monoid $M$, the {\em Weil restriction } of the functor $\ssp M$,  is defined by
\begin{equation}\label{resW}
\Ab \to \Se,\quad H\to\Hom_{\Mo}(M,\F_1[H]).
\end{equation}

The next proposition shows that the restriction to $\Ab$ of an $\Mo$-scheme is a direct sum of representable functors.

\begin{prop} \label{propdec} Let $\underline X$ be an $\Mo$-scheme and $X$ its geometric realization. Then the (restriction) functor
  $$
  \underline X: \Ab\to \Se,\quad H\mapsto \Hom(\Sp\F_1[H],X)=\underline X(\F_1[H])
  $$
  is the disjoint union
   \begin{equation}\label{disjunion}
  \underline X(\F_1[H])=\cup_{x\in X} \underline{X_x}(H)\,, \ \ \underline{X_x}(H)=\Hom_{\Ab} (\cO_x^\times,H).
  \end{equation}
   \end{prop}

\proof Let $\varphi \in  \Hom(\Sp\F_1[H],X)$. The unique point $\ffp\in \Sp\F_1[H]$ corresponds to the ideal $(0)$. Let
$\varphi(\ffp)=x\in X$ be its image; there is a corresponding map of the stalks
$$
\varphi^{\#}\,:\cO_{\varphi(p)}\to \cO_\ffp=\F_1[H].
$$
This homomorphism is local by hypothesis: this means that the inverse image of  $(0)$ by $\varphi^{\#}$ is
the maximal ideal of $\cO_{\varphi(p)}=\cO_x$. Therefore, the map $\varphi^{\#}$ is entirely determined by the group homomorphism $\rho\in\Hom_{\Ab}(\cO_x^\times,H)$ obtained as the restriction of $\varphi^{\#}$. Thus $\varphi \in  \Hom(\Sp\F_1[H],X)$ is entirely specified by a point $x\in X$ and a group homomorphism $\rho \in \Hom_{\Ab} (\cO_x^\times,H)$.
\endproof

In algebraic geometry, given an algebraic variety $X$ one knows that the degree of transcendence of the residue field $\kappa(x)$ of a point $x\in X$ measures the {\em dimension} of the closure $\overline{\{x\}}\subset X$. For an $\Mo$-scheme one has the following corresponding (local) notion

\begin{defn}\label{localdim} Let $X$ be a geometric $\Mo$-scheme and let $x\in X$ be a point. The {\em local dimension} $n(x)$ of $x$ is the rank of the abelian group $\cO_x^\times$.
\end{defn}

The local dimension determines a natural grading on the restriction of an $\Mo$-scheme $\underline X$ to $\Ab$ by assigning the degree $n(x)$ to the subset $\underline{X_x}(H)\subset \underline X(\F_1[H])$ in the decomposition \eqref{disjunion}.
\begin{prop} \label{propdecgrade}
The restriction to $\Ab$ of the following $\Mo$-schemes $\underline X$ coincides as a functor to $\Z_{\ge 0}$-graded sets with the functors defined in \cite{ak}
\begin{enumerate}
  \item Finite abelian groups $D$: $\underline X=\ssp \F_1[D]$
  \item Tori $\Gm$: $\underline X=\ssp \F_1[\Z]$
  \item Affine space $\cD^n$: $\underline X=\ssp \F_1[T_1,\ldots,T_n]$
  \item Projective line $\underline X=\spadu$
\end{enumerate}
\end{prop}

\proof 1) The space $\Sp \F_1[D]$ has a single point and the local dimension is $0$ which agrees with example 3.1 of \cite{ak}.

2) The space $\Sp \F_1[\Z]$ has a single point and the local dimension is $1$ which agrees with example 3.2 of \cite{ak}. This extends immediately to higher dimensional tori.

3)  Let $M=\F_1[T_1,\ldots,T_n]=\{0\}\cup\{T_1^{a_1}\cdots T_n^{a_n}|a_j\in\Z_{\ge 0}\}$.
 A prime ideal $\ffp$ of $M$ is of the form $\ffp = \cup_{i\in J}T_iM$, where $J$ is a subset of $\{1,\ldots,n\}$. One has $\mathcal O_\ffp^\times \simeq \Z^{J^c}$, generated by the $T_j$'s, with $j\notin J$. Thus the local dimension of $\Sp M$ at $\ffp$ is the cardinality of $J^c$ and this agrees
 with example 3.3  of \cite{ak}.

 4) The geometric realization $\spad$ is obtained by gluing two affine lines (\cf \cite{deit} and \S \ref{sectMo}) and consists of three points $\spad=\{0,u,\infty\}$
where  the local dimension is zero at $0$ and $\infty$ and is one at $u$. This agrees with example 3.4
of \cite{ak}.
\endproof

  \section{The  category $\Mr$ and $\F_1$-Schemes} \label{catmr}

 As we already remarked in \cite{ak} (\cf~\S 4), the definition of the (affine) variety over $\F_1$ for a Chevalley group is inclusive of the datum given by a covariant functor to the category $\Se$ of sets, fulfilling  much stronger properties than the ones required originally in \cite{Soule} (for affine varieties). The domain of such functor is a category which contains both the category of commutative rings and that of monoids (these two categories being linked by a pair of adjoint functors) and moreover in the definition of the variety one also requires the existence of a suitable natural transformation. In this section we develop the details of this construction following an idea we learnt from P.~Cartier. The introductory subsection \S \ref{adjfunctor} develops some generalities on the gluing process of two categories linked by a pair of adjoint functors. In \S \ref{extfunc}, we also treat in this generality the extension of functors.  The specific case of interest is covered in \S\ref{funfunctor} where we show that Chevalley groups are schemes over $\F_{1^2}$. Finally, in \S \ref{zetanoethersect} we extend the computation of zeta functions of \cite{deit2} (Theorem 1)  to our new setup which is no longer restricted to toric varieties (as it  covers in particular the case of Chevalley groups).\vspace{.05in}

 \subsection{Gluing two categories using adjoint functors}\label{adjfunctor}\hfill \vspace{.05in}

 We consider two categories $\cC$ and $\cC'$ and a pair of adjoint functors $\beta:\cC\to \cC'$
and $\beta^*:\cC'\to \cC$. Thus, one has a canonical identification
\begin{equation}\label{adjtrel0}
\Hom_{\cC'}(\beta(H), R)\stackrel{\Phi}{\cong} \Hom_{\cC}(H,\beta^*(R))\qquad\forall~H\in{\rm Obj}(\cC),~R\in{\rm Obj}(\cC').
\end{equation}
The naturality of  $\Phi$ is expressed by the commutativity of the following diagram where the vertical arrows are given by composition, $\forall f\in \Hom_{\cC}(G,H)$ and $\forall h\in \Hom_{\cC'}(R,S)$
\begin{gather} \raisetag{37pt} \,\hspace{60pt}
\xymatrix@C=25pt@R=25pt{
\Hom_{\cC'}(\beta(H),R) \ar[d]_{\Hom(\beta(f),h)} \ar[r]^{\Phi} &
\Hom_{\cC}(H,\beta^*(R)) \ar[d]^{\Hom(f,\beta^*(h))} \\
  \Hom_{\cC'}(\beta(G),S) \ar[r]_{\Phi}  & \Hom_{\cC}(G,\beta^*(S))
  } \label{diag}\hspace{100pt}
\end{gather}

We shall now define a category $\cC''=\cC\cup_{\beta,\beta^*} \cC'$ obtained by gluing $\cC$ and $\cC'$. The collection\footnote{It is not a set: we refer for details to the discussion contained in the preliminaries of \cite{demgab}} of objects of $\cC''$ is obtained as the {\em disjoint union} of the collection of objects of $\cC$ and $\cC'$. For $R\in{\rm Obj}(\cC')$ and $H\in{\rm Obj}(\cC)$, one sets $\Hom_{\cC''}(R,H)=\emptyset$. On the other hand, one defines
\begin{equation}\label{morC0}
    \Hom_{\cC''}(H,R)=\Hom_{\cC'}(\beta(H), R)\cong \Hom_{\cC}(H,\beta^*(R)).
\end{equation}
 The morphisms between objects contained in a same category are unchanged. The composition of morphisms in $\cC''$ is defined as follows. For $\phi \in \Hom_{\cC''}(H,R)$ and $\psi \in \Hom_{\cC}(H',H)$, one  defines $\phi\circ \psi\in \Hom_{\cC''}(H',R)$ as the composite
 \begin{equation}\label{defcomp1}
 \phi\circ \beta(\psi)\in \Hom_{\cC'}(\beta(H'), R)=\Hom_{\cC''}(H',R).
 \end{equation}
  Using the commutativity of the diagram  \eqref{diag}, one obtains
  \begin{equation}\label{defcomp1bis}
\Phi( \phi\circ \beta(\psi))=\Phi(\phi)\circ \psi\in \Hom_{\cC}(H',\beta^*( R)).
 \end{equation}
  Similarly, for $\theta\in \Hom_{\cC'}(R, R')$ one defines $\theta \circ \phi\in \Hom_{\cC''}(H,R')$ as the composite
  \begin{equation}\label{defcomp2}
 \theta \circ \phi\in \Hom_{\cC'}(\beta(H), R')=\Hom_{\cC''}(H,R')
 \end{equation}
 and using  again the commutativity of \eqref{diag} one obtains that
  \begin{equation}\label{defcomp2bis}
\Phi( \theta \circ \phi)=\beta^*(\theta)\circ\Phi(\phi)
\in \Hom_{\cC}(H,\beta^*( R')).
 \end{equation}

Moreover, one also defines specific morphisms $\alpha_H$ and $\alpha'_R$ as follows
\begin{equation}\label{alphas}
    \alpha_H={\rm id}_{\beta(H)}\in \Hom_{\cC'}(\beta(H), \beta(H))=\Hom_{\cC''}(H,\beta(H))
\end{equation}
\begin{equation}\label{alphasbis}
    \alpha'_R=\Phi^{-1}({\rm id}_{\beta^*(R)})\in \Phi^{-1}( \Hom_{\cC}(\beta^*(R), \beta^*(R)))=\Hom_{\cC''}(\beta^*(R),R).
\end{equation}
By construction one gets
\begin{equation}\label{expressbis}
    \Hom_{\cC''}(H,R)=\{g\circ\alpha_H\,|\, g\in \Hom_{\cC'}(\beta(H),R)\}
\end{equation}
and for any morphism $\rho\in \Hom_{\cC}(H,K)$ the following equation holds
\begin{equation}\label{exprrel}
\alpha_K\circ \rho=\beta(\rho)\circ \alpha_H.
\end{equation}
Similarly, it also turns out that
\begin{equation}\label{express}
    \Hom_{\cC''}(H,R)=\{\alpha'_R\circ f\,|\, f\in \Hom_{\cC}(H,\beta^*(R))\}
\end{equation}
and the associated equalities hold
\begin{equation}\label{exprrelbis}
\alpha'_S\circ \beta^*(\rho)= \rho\circ \alpha'_R \qquad\forall \rho\in \Hom_{\cC'}(R,S)
\end{equation}
\begin{equation}\label{relate}
    g\circ\alpha_H=\alpha'_R\circ \Phi(g)\qquad\forall g\in \Hom_{\cC'}(\beta(H),R).
\end{equation}

\vspace{.05in}

\begin{prop}\label{catplus0}  $\cC''=\cC\cup_{\beta,\beta^*} \cC'$ is a category which contains $\cC$ and $\cC'$ as full subcategories. Moreover, for any object $H$ of $\cC$ and $R$ of $\cC'$, one has
$$
 \Hom_{\cC''}(H,R)=\Hom_{\cC'}(\beta(H), R)\cong \Hom_{\cC}(H,\beta^*(R)).
$$
\end{prop}

\proof  One needs to check that the composition $\circ''$ of morphisms is associative in $\cC''$ \ie that
$h \circ'' (g \circ'' f) = (h \circ'' g) \circ'' f$. The only relevant case to check is when the image of $f$ is an object $H$ of $\cC$ and the image of $g$ is an object $R$ of $\cC'$.
 Then $f(G)=H$, with $G$ an object of $\cC$ and $h(R)=S$ an object of $\cC'$. One has $g\in \Hom_{\cC'}(\beta(H), R)$ and
 $$
 g \circ'' f=g\circ \beta(f)\,, \ h \circ'' (g \circ'' f)=h\circ(g\circ \beta(f))=
 (h\circ g)\circ \beta(f)=(h \circ'' g) \circ'' f
 $$
\endproof

\subsection{Extension of functors.}\label{extfunc}\hfill \vspace{.05in}

 We keep the notations introduced in \S\ref{adjfunctor} and let $(\beta,\beta^*)$ be a pair of adjoint functors linking $\cC$ and $\cC'$, \ie $\beta: \cC\to \cC'$ and $\beta^*: \cC'\to \cC$ and the isomorphism \eqref{adjtrel0} holds.  Let $\cF: \cC \to \cT$ and $\cF': \cC' \to \cT$ be covariant functors   to the same category $\cT$. It is straightforward routine to verify that the assignment of a natural transformation $\cF\to\cF'\circ \beta$ is equivalent to giving a natural transformation $\cF\circ \beta^*\to\cF'$. By implementing in this setup the category $\cC''=\cC\cup_{\beta,\beta^*} \cC'$ defined in \S\ref{adjfunctor}, one obtains the following more precise result

\begin{prop}\label{catplus}  1)  With the above notation, let $\cF'':\cC''\to\cT$ denote a covariant functor. Then the  assignment $H\mapsto\cF''(\alpha_H)$ defines a
natural transformation $\cF''|_{\cC}\to\cF''|_{\cC'}\circ \beta$
 and analogously the assignment $R\mapsto\cF''(\alpha'_R)$ defines a natural transformation $\cF''|_{\cC}\circ \beta^*\to\cF''|_{\cC'}$.\vspace{.05in}

 2) Let $\cF: \cC \to \cT$ and $\cF': \cC' \to \cT$ be covariant functors. Then \vspace{.05in}

 a) Given a  natural transformation  $\cF\to\cF'\circ
\beta$, there exists a unique covariant functor $\cF''$  which extends $\cF$,
$\cF'$ and
 agrees with the natural transformation on the morphisms
$\alpha_H$.

b) Given a   natural transformation $\cF\circ \beta^*\to\cF'$, there exists a unique covariant functor $\cF''$  which extends $\cF$,
$\cF'$ and
 agrees with the natural transformation on the morphisms $\alpha'_R$.
\end{prop}

\proof $1)$ follows from \eqref{exprrel} and \eqref{exprrelbis}.

$2)\, a)$ A natural transformation $\cF\to\cF'\circ \beta$ determines, by
 \eqref{expressbis}, the extension   from $\cC\cup \cC'$ to $\cC''=\cC\cup_{\beta,\beta^*} \cC'$.

$2)\, b)$  The proof is similar to the proof of $2 a)$.
\endproof

Let us assume now  that we are given a functor $\cF: \cC \to \cT$, where  $\cT=\Se$. In the following we shall investigate under which  conditions $\cF$  admits an extension to
 $\cC''=\cC\cup_{\beta,\beta^*} \cC'$.

 Notice first that if $\cF$ is representable, then it admits a unique representable extension to
$\cC''$. Indeed, the representability of $\cC$ amounts to the existence of an object $G$ in $\cC$ such that
$$
\cF(H)=\Hom_{\cC}(G,H)\qquad\forall H\in{\rm Obj}(\cC).
$$
If the extension of $\cF$ is represented by an object of $\cC''$ then this object must necessarily belong to $\cC$, since by definition of $\cC''$ there is no morphism of $\cC''$ from an object  of $\cC'$ to an object of $\cC$. Moreover, by restriction to $\cC$ one gets the uniqueness of the extension. Thus, any extension of $\cF$ to $\cC''$ which is represented by an object of this latter category is unique (as a representable functor) and is given on $\cC'$ by
\begin{equation}\label{reprepun}
\cF'(R)=\Hom_{\cC'}(\beta(G),R).
\end{equation}
The natural transformation $\cF\to\cF'\circ \beta$ is simply given by the restriction of the functor $\beta$
$$
\beta:\Hom_{\cC}(G,H)\to \Hom_{\cC'}(\beta(G),\beta(H)).
$$
Similarly, the natural transformation $\cF\circ \beta^*\to\cF'$ is defined by the identity map
$$
\Hom_{\cC}(G,\beta^*(R))\to \cF'(R)=\Hom_{\cC'}(\beta(G),R)\cong \Hom_{\cC}(G,\beta^*(R)).
$$
Thus, the following result holds

\begin{prop}\label{catext}  Let $\cF: \cC \to \Se$ be a representable functor.\vspace{.05in}

1) There exists a unique extension $\tilde\cF$ of $\cF$ to $\cC''=\cC\cup_{\beta,\beta^*} \cC'$ as a representable functor.\vspace{.05in}

2) Let $\cG$ be any extension of $\cF$ to $\cC''=\cC\cup_{\beta,\beta^*} \cC'$. Then, there exists a unique morphism of functors from $\tilde\cF$ to $\cG$ which restricts to the identity on $\cC$.
\end{prop}

\proof For $1)$, we notice that the object of $\cC$ representing $\cF$ is unique up to isomorphism and it represents
$\tilde\cF$.

The proof of $2)$ follows from the following facts. If $A\in{\rm Obj}(\cC)$ represents $\tilde\cF$, by applying Yoneda's lemma (\ie ${\rm Nat}(\tilde\cF,\cG) \simeq \cG(A)$) we know that there exists a uniquely determined  natural transformation $\phi: \tilde\cF\to\cG$ associated to any object $\rho$ of $\cG(A)$: the pair $(\phi,\rho)$  being linked by the formulas $\rho=\phi(A)({\rm id}_A)\in \cG(A)$ and $\phi = \rho^\sharp$ ($\rho^\sharp(R)(\beta) = \cG(\beta)(\rho)$).  But since the restriction of $\phi$ to $\cC$ is the identity map from $\cF(A)$ to $\cG(A)=\cF(A)$, one obtains the required uniqueness.
\endproof

The following corollary shows that even though the extension  $\tilde\cF$ of $\cF$ to $\cC''=\cC\cup_{\beta,\beta^*} \cC'$ is not unique it is universal.

\begin{cor}\label{coruniv} Let $\cF: \cC \to \Se$ be a representable functor and let $\cF': \cC' \to \Se$ be defined as in \eqref{reprepun}. Let $\cG'$ be a functor from $\cC'$ to $\Se$ and let $\phi: \cF \to \cG'\circ \beta$ be a natural transformation. Then there exists a unique morphism of functors $\psi:\cF'\to \cG'$ such that
\begin{equation}\label{basic1}
    \phi_H=\psi_{\beta(H)}\circ\tilde\cF(\alpha_H)\qquad\forall H\in{\rm Obj}(\cC).
\end{equation}
\end{cor}

\proof Given $\phi$  and $\cG'$, there exists by Proposition \ref{catplus} a unique extension $\cG$ of $\cF$ to $\cC''=\cC\cup_{\beta,\beta^*} \cC'$ which restricts to $\cG'$ on $\cC'$ and is such that
$$
\phi_H=\cG(\alpha_H)\qquad\forall H\in{\rm Obj}(\cC).
$$
A morphism of functors from $\tilde\cF$ to $\cG$ extending the identity on $\cC$ is entirely specified by its restriction to $\cC'$ which is a morphism of functors $\psi$ from $\cF'$ to $\cG'$ and it must be compatible with the morphisms $\alpha_H$. This compatibility is given by
\eqref{basic1}. Thus the existence and uniqueness of $\psi$ follows from Proposition \ref{catext}.
\endproof

The next proposition states a similar, but simpler, result for extensions of functors from
$\cC'$ to the larger category $\cC''=\cC\cup_{\beta,\beta^*} \cC'$.

\begin{prop}\label{catextbis}  Let $\cF':\cC' \to \Se$ be a  functor.\vspace{.05in}

1) There exists a unique extension $\tilde\cF$ of $\cF'$ to $\cC''=\cC\cup_{\beta,\beta^*} \cC'$ given by $\cF'\circ \beta$ on $\cC$ and such that $\tilde\cF(\alpha_H)={\rm id}_{\beta(H)}$ for all objects $H$ of $\cC$.\vspace{.05in}

2) Let $\cG$ be any extension of $\cF'$ to $\cC''=\cC\cup_{\beta,\beta^*} \cC'$, then there exists a unique morphism of functors from $\cG$ to $\tilde\cF$   which is the identity on $\cC'$.
\end{prop}

\proof The first statement follows from $2)$ of Proposition \ref{catplus} by using the identity as a natural transformation. Similarly, for the second statement, again $2)$ of Proposition \ref{catplus} determines a unique morphism of functors $\phi$ given by $\phi(H)=\cG(\alpha_H)$ and obtained from the restriction of $\cG$ to $\cC$ to $\cG\circ \beta=\cF'\circ \beta$. We extend $\phi$ as the identity on $\cC'$. The  commutative diagram
\begin{gather} \raisetag{37pt} \,\hspace{100pt}
\xymatrix@C=25pt@R=25pt{
\cG(H) \ar[d]^{\phi_H} \ar[r]^{\cG(\alpha_H)} &
\cG(\beta(H)) \ar[d] \\
  \tilde \cF(H)=\cF'(\beta(H)) \ar[r]  & \tilde \cF(\beta(H))
  } \label{diagbis}\hspace{100pt}
\end{gather}
shows that one obtains a morphism of functors from $\cG$ to $\tilde\cF$. The same diagram also gives the uniqueness of $\phi$.
\endproof\vspace{.05in}

 \subsection{$\F_1$-Schemes and Chevalley groups.}\label{funfunctor}\hfill \vspace{.05in}

 We now apply the  construction described in \S\S\ref{adjfunctor}, \ref{extfunc} to the pair of adjoint covariant functors $\beta$ and $\beta^*$ which are defined as follows. The functor
\begin{equation}\label{beta}
 \beta: \Mo \to \An\,, \quad   M\mapsto \beta(M)=\Z[M]
\end{equation}
associates to a monoid $M$ the convolution ring $\Z[M]$ (the $0$ element of $M$ is sent to $0$). The adjoint functor $\beta^*$
\begin{equation}\label{units}
 \beta^*: \An \to \Mo\quad   R\mapsto \beta^*(R)= R
\end{equation}
 associates to a ring $R$ the ring itself viewed as a multiplicative monoid. The adjunction relation means that
\begin{equation}\label{adjtrel}
\Hom_\An(\beta(M), R)\cong \Hom_\Mo(M,\beta^*(R)).
\end{equation}
Now, we apply Proposition \ref{catplus0} to construct the category $\Mr=\An\cup_{\beta,\beta^*} \Mo$.
 Thus, one obtains for every object $R$ of $\An$, a morphism
\begin{equation}\label{added}
    \alpha'_R\in \Hom_\Mr(\beta^*(R),R)
\end{equation}
and the following relation between the morphisms of $\Mr$:
\begin{equation}\label{glue}
  f\circ  \alpha'_R=\alpha'_S\circ\beta^*(f)\qqq f\in \Hom_\An(R,S).
\end{equation}
Similarly, for every  monoid $M$ one has a morphism
\begin{equation}\label{addedm}
    \alpha_M\in  \Hom_\Mr(M,\beta(M))
\end{equation}
together with the relation
\begin{equation}\label{gluebis}
  \beta(f)\circ  \alpha_M=\alpha_N\circ f\qqq f\in \Hom_\Mo(M,N).
\end{equation}

\begin{defn}\label{defnfunfunc} An $\F_1$-functor is a covariant functor from the category
$\Mr=\An\cup_{\beta,\beta^*} \Mo$ to the category of sets.
\end{defn}

Then, it follows from Proposition~\ref{catplus} that giving an $\F_1$-functor $\mathcal X: \Mr \to \mathfrak{Sets}$ is equivalent to the assignment of the following data:\vspace{.05in}

$\bullet$~An $\Mo$-functor $\underline X$.\vspace{.05in}

$\bullet$~A $\Z$-functor $X_\Z$.\vspace{.05in}

$\bullet$~A natural transformation $e:\underline X\to X_\Z\circ \beta$.\vspace{.05in}

The third condition can be equivalently replaced by the assignment of a natural transformation $\underline X\circ \beta^*\to X_\Z$.\vspace{.05in}

Now that we have at our disposal the category $\Mr$ obtained by gluing $\Mo$ and $\An$ we introduce our notion of an $\F_1$-scheme.

\begin{defn}\label{defnfonesch} An $\F_1$-scheme is an $\F_1$-functor $\mathcal X: \Mr\to \mathfrak{Sets}$, such that:\vspace{.05in}

$\bullet$~The restriction $X_\Z$ of $\mathcal X$ to $\An$ is a $\Z$-scheme.\vspace{.05in}

$\bullet$~The restriction $\underline X$ of $\mathcal X$ to $\Mo$ is an $\Mo$-scheme.\vspace{.05in}

$\bullet$~The natural transformation $e:\underline X\circ \beta^*\to X_\Z$ associated to a field is a bijection (of sets).
\end{defn}

Morphisms of $\F_1$-schemes are natural transformations of the corresponding functors.\vspace{.05in}

Notice that even though the restriction $\underline X$ of $\mathcal X$ to $\Mo$  is an $\Mo$-scheme, the composite $\underline X\circ\beta^*$ is not in general a $\Z$-scheme since it is not a local $\Z$-functor. As an example, one may consider the case of $\mathcal X=\P^1$: here, the above composition of functors determines only a smaller portion of the projective line as a $\Z$-scheme. However, one can associate to $\underline X\circ\beta^*$ a {\em unique} $\Z$-scheme $Y$ that is defined by assigning to a ring $R$ the set $Y(R)$ of solutions of \eqref{locf}, using $\underline X\circ\beta^*$ and an arbitrary partition of unity in $R$. Proposition \ref{catext} describes a {\em canonical} morphism of $\Z$-schemes from $Y$ to $X_\Z$. However, this morphism is not in general an isomorphism: we refer, for instance, to the case of Chevalley groups described in \cite{ak}. Proposition~\ref{catextbis} describes the natural morphism from the $\Mo$-scheme $\underline X$ to the ``gadget'' (\cf~\cite{ak} Definition~2.5) associated to the $\Z$-scheme $X_\Z$.\vspace{.05in}

Definition \ref{defnfonesch} admits variants corresponding to the extensions $\F_{1^n}$. For Chevalley schemes we shall  need to consider only the case $n=2$. One defines the category $\Mo^{(2)}$ of pairs made by an object $M$ of $\Mo$ and an element $\epsilon\in M$ of square one. One defines the functors $\beta: \Mo^{(2)}\to \An$ as $\beta(M,\epsilon)=\Z[M,\epsilon]$
 \begin{equation}\label{constadj}
\beta(M,\epsilon)=\Z[M,\epsilon]=\Z[M]/J\,, \ \ J=(1+\epsilon)\Z[M]
\end{equation}
and $\beta^*: \An\to \Mo^{(2)}$ as $\beta^*(A)=(A,-1)$, where $(A,-1)$ is the object of $\Mo^{(2)}$ given by the ring $A$ viewed as a
(multiplicative) monoid and the element $-1\in A$. Then, the following adjunction relation holds for any commutative ring $A$
\begin{equation}\label{adjadj}
    \Hom_\An(\beta(M,\epsilon),A)\cong\Hom_{\Mo^{(2)}}((M,\epsilon),\beta^*(A)).
\end{equation}

Finally, we have the following result

\begin{thm} \label{mainthmbis}  The  scheme ${\mathfrak
  G}$ over $\Z$ associated to a Chevalley group $G$ extends to a scheme $\cG$ over $\F_{1^2}$.
\end{thm}
\proof The proof follows from \cite{ak} Theorem~4.1, where we showed that

$\bullet$~The construction  of the functor $\underline G$ extends from the category $\mathcal F_{ab}^{(2)}$ of pairs $(D,\epsilon)$ of a finite abelian group and an element of order two to
the category $\Mo^{(2)}$.\vspace{.05in}

$\bullet$~The construction  of the natural transformation $e_G$ extends from $\mathcal F_{ab}^{(2)}$ to the category $\Mo^{(2)}$ and associates to any   $A\in{\rm Obj}(\An)$  a map

\begin{equation}\label{mapA}
  e_{G,A} :\;\underline G(A,-1) \to \mathfrak G(A)
  \end{equation}

$\bullet$~ When $A$ is a field the map $ e_{
  G,A}$ is a bijection.\endproof\vspace{.05in}

The map $ e_{G,A}$ of \eqref{mapA} is constructed in \cite{ak}, \cf~proof of Theorem~4.1, it yields the natural transformation from $\underline G\circ \beta^*$ to $\mathfrak    G$ (and using \eqref{relate} the corresponding natural transformation from $\underline G$ to $\mathfrak  G\circ \beta$). For $H\in {\rm Obj}(\Ab)$, this natural transformation is compatible with the group structure on the subset $\underline G^{(\ell)}(H)$. We  recall that $\underline G^{(\ell)}(H)\subset \underline G(H)$ is the part of lowest degree ($\ell={\rm rk}~\mathfrak G$) in the grading of $\underline G$ (\cf Definition \ref{localdim}).

\subsection{Zeta function of Noetherian $\F_1$-schemes}\hfill\label{zetanoethersect}\vspace{.05in}

 We recall that a congruence on a monoid $M$ is an equivalence relation which is compatible with the semigroup operation. A  monoid is Noetherian when any strictly increasing sequence of congruences is finite (\cf~\cite{Gilmer} page 30).
The following conditions on a monoid $M$ are equivalent (\cf~\cite{Gilmer}: Theorems~7.7, 7.8 and 5.10)\vspace{.05in}

-~$M$ is Noetherian.\vspace{.05in}

-~$M$ is finitely generated.\vspace{.05in}

-~$\Z[M]$ is a Noetherian ring.\vspace{.05in}

It is proven in \cite{Gilmer} (Theorem 5.1) that if $M$ is a Noetherian monoid, then for any prime ideal $\ffp\subset M$ the localized monoid $M_\ffp$ is also Noetherian (the semigroup $\ffp^c$ is finitely generated and thus $M_\ffp$ is also finitely generated). The same theorem also shows that the abelian group $(M_\ffp)^\times$ is finitely generated.

 \begin{defn}\label{noetherdefn}
 An $\Mo$-scheme is Noetherian if it admits a finite open cover by representable subfunctors $\ssp(M)$,  with $M$ Noetherian monoids.

 An $\F_1$-scheme is Noetherian if the associated $\Mo$ and $\Z$-schemes are Noetherian.
\end{defn}

A geometric $\Mo$-scheme $X$ is said to be torsion free if the groups $\cO_x^\times$ of invertible elements of the monoids $\cO_x$, for $x\in X$ are torsion free. The following result is related to Theorem 1 of \cite{deit2}, but it applies to a wider class of varieties (\ie non necessarily toric)

\begin{thm}\label{dthmfonesch}
Let $\mathcal X$ be a Noetherian $\F_1$-scheme and let $X$ be the geometric realization of its restriction $\underline X$ to $\Mo$. Then, if $X$ is torsion free,
\begin{enumerate}
  \item There exists a polynomial $N(x+1)$ with positive integral coefficients such that
  $$
  \# \,\underline X(\F_{1^n})=N(n+1)\qqq n\in \N.
  $$
  \item For each finite field $\F_q$ the cardinality of the set of points of the $\Z$-scheme
  $X_\Z$ which are rational over $\F_q$ is equal to $N(q)$.
  \item The zeta function of $\mathcal X$  has the following description
  \begin{equation}\label{prodfund}
    \zeta_\mathcal X(s)=\prod_{x\in X} \frac{1}{\left(1-\frac 1 s\right)^{\otimes^{n(x)}}}
\end{equation}
where $\otimes$ denotes Kurokawa's tensor product and  $n(x)$ is the local dimension of $X$ at the point $x$.
\end{enumerate}
\end{thm}

In \eqref{prodfund}, we use the convention that when $n(x)=0$ the expression
$$
\left( 1-\frac 1 s\right)^{\otimes^{n(x)}}=s.
$$
We refer to \cite{Ku} and \cite{Manin} for the details of the definition of Kurokawa's tensor products and zeta functions. When $\mathcal X$ is the projective line $\spad$, one has three points in $X$: two are of dimension zero and one is of dimension one. Thus  \eqref{prodfund} gives
$$
\prod_{x\in X} \frac{1}{\left(1-\frac 1 s\right)^{\otimes^{n(x)}}}=\frac{1}{s^2}\frac{1}{1-\frac 1 s}=\frac{1}{s(s-1)}.
$$
Formula \eqref{prodfund} continues to hold even in the presence of torsion  and in that case it corresponds to the treatment of torsion given in \cite{deit2}.

\proof $(1)$ By construction,  $\underline X(\F_{1^n})$ is the set obtained by evaluating the restriction of $\underline X$ (from the subcategory $\Mo$ of $ \Mr$) on $\Ab$, at the cyclic group $H=\Z/n\Z$. By Proposition \ref{propdec} one has
 \begin{equation}\label{formpropdec}
 \underline X(H)=\cup \underline X_x(H)\,, \ \ \underline X_x(H)=\Hom_{\Ab} (\cO_x^\times,H).
 \end{equation}
Since $X$ is Noetherian it is a finite topological space  and for each $x\in X$ the abelian group $\cO_x^\times$ is finitely generated and torsion free by hypothesis. The rank of $\cO_x^\times$ is $n(x)$ and thus the set
$\underline X_x(H)=\Hom_{\Ab} (\cO_x^\times,H)$ has cardinality $n^{n(x)}$. It follows that
$$
P(y)=\sum_{x\in X} y^{n(x)}.
$$
is a polynomial function in the indeterminate $y$ with positive integral coefficients which is related to the counting function of $\mathcal X$ by the equation $N(x+1)=P(x)$. Then $(1)$ follows.

$(2)$ follows from $(1)$ and the fact that the natural transformation $e$ (which is part of the set of data describing an $\F_1$-scheme, \cf~Definition~\ref{defnfonesch}) evaluated at any field is a bijection. In the case of a finite field $\F_q$, the corresponding monoid is $\F_1[H]$ for the cyclic group $H=\Z/n\Z$ of order $n=q-1$.

To prove $(3)$, we start by computing explicitly the Kurokawa's tensor product (for $n>0$) as follows
\begin{equation}\label{kuro}
\left( 1-\frac 1 s\right)^{\otimes^{n}}=\prod_{j\, {\rm even}} (s-n+j)^{n\choose j}/\prod_{j\, {\rm odd}} (s-n+j)^{n\choose j}.
\end{equation}
The above equality is a straightforward consequence of the definition of the Kurokawa's tensor product since the divisor of zeros of
$1-\frac 1 s$ is $\{1\}-\{0\}$ and its $n$-th power is given by the binomial formula
$$
(\{1\}-\{0\})^n=\sum (-1)^k{n\choose k}\{n-k\}.
$$
To obtain \eqref{prodfund} we use \eqref{formpropdec} to express $\zeta_\mathcal X(s)$ as a product over the points of $X$ and then we apply the well-known fact
 \begin{equation}\label{polynomial}
N(x)=\sum_0^d a_kx^k ~\Longrightarrow~ \zeta_\mathcal X(s)=\prod_0^d(s-k)^{-a_k}
\end{equation}
and \eqref{kuro} to show that the zeta function for the  polynomial $(q-1)^n$ is the inverse of $\left( 1-\frac 1 s\right)^{\otimes^{n}}$.
\endproof

\section{The projective ad\`ele class space}\label{projadel}

In this section we develop an application of the functorial approach of the theory of $\Mo$-schemes described in \S \ref{sectMo} to explain, at a conceptual level, the spectral realization of zeros of $L$-functions for an arbitrary global field $\K$. \vspace{.05in}

\subsection{Vanishing result for $\Mo$-schemes}\label{vanmo}
\hfill\vspace{.05in}

In this subsection, we shall first briefly review some standard facts on sheaf cohomology and then we will show that for sheaves of abelian groups over (the geometric realization of) an $\Mo$-scheme, sheaf cohomology and \v{C}ech cohomology agree. \vspace{.05in}

 Given a topological space $X$, we denote by $\Ab(X)$ the category of sheaves of abelian groups on $X$. It is a well-known fact that $\Ab(X)$ is an abelian category with enough injectives (\cf \cite{tohoku}, Propositions 3.1.1 and  3.1.2).  For any open subset $U\subset X$ the functor
 $$
 \Gamma(U,\cdot):\Ab(X) \longrightarrow \Ab,\quad F\mapsto \Gamma(U,F)=F(U)
 $$
  describes the space of sections of $F$ on $U$ and it is left-exact. Its derived functor defines the sheaf cohomology  $H(U,F)$. Moreover, for any point $x\in X$ the functor `stalk of $F$ at $x$'
\begin{equation}\label{exactstalk}
 \Ab(X) \longrightarrow \Ab,\quad F\mapsto \varinjlim_{x\in U}\Gamma(U,F)=:F_x
\end{equation}
is exact.

\begin{prop}\label{cechprop} Let $X$ be the geometric realization of an $\Mo$-scheme, then the following results hold\vspace{.05in}

1) For any open affine set $U\subset X$
\begin{equation}\label{vanish}
    H^p(U,F)=0\qquad\quad \forall p>0,\quad \forall~F\in {\rm Obj}\;\Ab(X).
\end{equation}
2) Let $\cU=\{U_j\}_{j\in J}$ be an open cover of $X$ such that all finite intersections $\bigcap_{j_k} U_{j_k}$
are affine, then for any sheaf $F$ of abelian groups  on $X$, one has
\begin{equation}\label{cechequal}
    H^p(X,F)=\check{H}^p(\cU,F)\qquad\forall p\geq 0
\end{equation}
where the cohomology on the right hand side is the \v{C}ech cohomology relative to the covering $\cU$.

3) Let $Y=U^c$ be the complement of an affine open set in $X$. Then for any sheaf $F$ of abelian groups on $X$ one has the exact sequence
\begin{eqnarray}
  0 &\to & H_Y^0(X,F)\to H^0(X,F)\to H^0(U,F|_{U})\label{longexact1} \\
   &\to & H_Y^1(X,F)\to H^1(X,F)\to 0 \nonumber
\end{eqnarray}
where $H^*_Y(X,F)$ denotes the cohomology with support on $Y$.
\end{prop}

\proof $1)$ Let $U=\Sp M$, where $M$ is a monoid in $\Mo$. Then any open set $V\subset U$ which contains the closed point $\ffp=(M^\times)^c$ of $U$ coincides with $U$ (\cf~\S 3.2 Lemma~3.3). Thus the stalk $F_\ffp$ is equal to $\Gamma(U,F)=F(U)$ and the result follows from the exactness of the functor `stalk at $\ffp$' \eqref{exactstalk}.

$2)$ follows from $1)$ in view of the equality of $H^p(X,F)$ with the \v{C}ech cohomology relative to the covering $\cU$, under the assumption that for  all finite intersections $V=\bigcap_{j_k} U_{j_k}$ of opens in $\cU$ one has $H^p(V,F)=0\qqq p>0$ (\cf \cite{Hart}, Exercice III, 4.11).

$3)$ Following \cite{Hart} III 2.3, one has for  any sheaf $F$ of abelian groups on $X$, a long exact sequence of cohomology groups
\begin{eqnarray}\label{long}
  0 &\to & H_Y^0(X,F)\to H^0(X,F)\to H^0(U,F|_U) \\
   &\to & H_Y^1(X,F)\to H^1(X,F)\to H^1(U,F|_U) \nonumber\\
   &\to & H_Y^2(X,F)\to \cdots \nonumber
\end{eqnarray}
where $F|_U$ denotes the restriction of the sheaf $F$ on the open set $U$. Thus 3) follows from \eqref{long} and the vanishing of $H^1(U,F|_U)$ shown in 1).
\endproof\vspace{.05in}

\subsection{The map to the base}\hfill\medskip\vspace{.05in}

We recall that an $\Mo$-scheme is (in particular) a covariant functor $\Mo\to\Se$, $M\to\underline X(M)$ and that there exists a {\em unique} geometric space $X$ associated to  $\underline X$ and satisfying the property that $\underline X(M) = \text{Hom}(\Sp(M),X)$.  For each monoid $M$, we let $\ffp_M$ be its maximal ideal: this is the complement of $M^\times$ in $M$.

\begin{prop} \label{coverspecbis}
Let
 $\underline X$ be an $\Mo$-scheme and let $X$ be its geometric realization.\vspace{.05in}

1) For any monoid $M$ there exists a canonical map of sets
\begin{equation}\label{pim}
    \pi_M\,:\, \underline X(M)\to X
\end{equation}
such that
\begin{equation}\label{pim1}
\pi_M(\phi)=\phi(\ffp_M),\qquad \phi\in \Hom(\Sp(M),X).
\end{equation}

2) Let $U$ be an open subset of $X$ and let $\underline U$ be the associated open subfunctor of $\underline X$, then
 \begin{equation}\label{pim2}
    \underline U(M)=\pi_M^{-1}(U)\subset \underline X(M).
 \end{equation}
  \end{prop}

  \proof $1)$ is a definition of the map $\pi_M$.

  $2)$ holds since for $\phi\in \underline X(M)=\Hom(\Sp(M),X)$ one has
  $$
  \phi(\ffp_M)\in U\iff \phi^{-1}(U)=\Sp(M).
  $$
  \endproof\vspace{.05in}

\subsection{The monoid $M=\A_\K/\K^\times$ of ad\`ele classes}  \hfill\vspace{.05in}

Let $\K$ be a  global field.
The id\`ele  class group $C_\K$ is the group $M^\times$ of invertible elements of the monoid
\begin{equation}\label{monoideclass}
 M=\A_\K/\K^\times\,, \qquad \K^\times=\GL_1(\K).
\end{equation}
We consider the $\Mo$-functor $\spadu$ associated to the projective line  of Example~\ref{projex}.
The geometric realization of $\spadu$ (\cf \cite{deit} and \S \ref{sectMo} above) is
the finite topological space $\spad$  whose underlying set is made by three points $\spad=\{0,u,\infty\}$
\begin{equation}\label{projspace1}
\overline{\{0\}}=\{0\},~ \ \overline{\{u\}}=\spad, ~ \overline{\{\infty\}}=\{\infty\}
\end{equation}
and whose topology is described by the following three open sets
\begin{equation}\label{projspace2}
U_+=\spad\backslash\{\infty\}\,, \ \ U_-=\spad\backslash\{ 0\}\,, \ \ U=U_+\cap U_-.
\end{equation}
We apply this construction to $M=\A_\K/\K^\times$ and obtain the {\em projective ad\`ele class space} $\spadu(M)$.
The canonical  projection  \eqref{pim} is given by
\begin{equation}\label{projback}
    \pi_M\,:\,   \spadu(M)=M\cup_{M^\times}M\to \spad.
\end{equation}
 It maps an element of $M^\times=C_\K$ to $u\in \spad$ and the complement of $M^\times=C_\K$  to $0$ or $\infty$ accordingly to the two copies of $M\backslash M^\times$ inside $\spadu(M)$.\vspace{.05in}

\subsection{The space of functions on the projective ad\`ele class space}\label{functdefn}\hfill\vspace{.05in}

 To define a natural space  $\cS(M) $ of functions on the quotient space $M=\A_\K/\K^\times$ of ad\`ele classes we consider the Bruhat-Schwartz space $\cS(\A_\K)$ of the locally compact abelian group $\A_\K$ and the space of its coinvariants under the action of $\K^\times$ (\cf\cite{Co-zeta}). More precisely we start with the exact sequence associated to the kernel
 of the $\K^\times$-invariant linear forms $\epsilon(f)=(f(0),\int_{\A_\K} f(x)dx)\in \C\oplus \C[1]$,
\begin{equation}\label{fonction1}
0\to \cS(\A_\K)_0 \to \cS(\A_\K)\stackrel{\epsilon}{\to} \C\oplus \C[1]\to 0\,.
\end{equation}
 Then one lets
\begin{equation}\label{fonction3}
\cS(M) = \cS_0(M)\oplus \C\oplus \C[1]\,, \ \  \cS_0(M)=\cS(\A_\K)_0 /\overline{\{f-f_q\}}
\end{equation}
where $\overline{\{f-f_q\}}$ denotes  the closure of the subspace of $\cS(\A_\K)_0$ generated by the differences $f-f_q$, with $q\in\K^\times$.\vspace{.05in}

We now introduce the functions on  the projective ad\`ele class space $\spadu(M)$.
We start by defining the sheaf $\Omm$ on $\spad$ which is {\em uniquely} determined by the following spaces of sections and restriction maps
\begin{eqnarray}\label{sheaf1}
   \Gamma(U_+,\Omm) &=& \cS(M)\nonumber \\
   \Gamma(U_-,\Omm) &=& \cS(M) \nonumber \\
   \Gamma(U_+\cap U_-,\Omm) &=& \strong(C_\K) \nonumber
\end{eqnarray}
where, for a number field $\K$, $\strong(C_\K)$ is defined as follows
\begin{equation}\label{SCK}
\strong(C_\K)=\,\cap_{\beta \in \R}\,\mu^\beta \cS(C_\K)=\{f\in \cS(C_\K)\,|\,\mu^\beta f\in \cS(C_\K)\qqq \beta \in \R\}.
\end{equation}
Here $\mu \in C(C_\K)$ denotes the module $\mu:C_\K\to\R_+^\times$, $\mu^\beta(g)=\mu(g)^\beta$ and $ \cS(C_\K)$ is the Bruhat-Schwartz space over $C_\K$ (\cf \cite{CMbook}, \cite{Meyer}).
The restriction maps to $U=U_+\cap U_-$ vanish on the component $\C\oplus \C[1]$ while on the other components they are defined as follows
\begin{eqnarray}\label{restmaps}
  (\Res\, f)(g) &=&  \sum_{q\in \K^\times}f(q g), \qquad f\in \cS_0(M)\subset\Gamma(U_+,\Omm)\nonumber  \\
  (\Res\, h)(g) &=& |g|^{-1}\sum_{q\in \K^\times}h(q g^{-1}),
  \qquad h\in \cS_0(M)\subset\Gamma(U_-,\Omm).
\end{eqnarray}\vspace{.05in}

We also introduce the notation:
\begin{equation}\label{Smap}
    \Sigma(f)(g)=\sum_{q\in \K^\times}f(qg).
\end{equation}\vspace{.05in}

\subsection{$H^0(\spad,\Omm)$ and the graph of the Fourier transform}\hfill\vspace{.05in}

The  \v{C}ech complex of the covering $\cU=\{U_\pm\}$ of $\P^1_{\F_1}$ has two terms
\begin{eqnarray}
  C^0 &=&  \Gamma(U_+,\Omm)\times  \Gamma(U_-,\Omm) \nonumber \\
  C^1 &=& \Gamma(U_+\cap U_-,\Omm). \nonumber
\end{eqnarray}
The coboundary $\partial:C^0\to C^1$ is given by (\cf~\eqref{Smap})
\begin{equation}\label{bord}
    \partial(f,h)(g)=\Sigma(f)(g)-|g|^{-1}\Sigma(h)(g^{-1}).
\end{equation}

\begin{lem}\label{lemsurject}
The kernel of $\partial$ is the graph of the Fourier transform \ie \begin{equation}\label{hzero}
    H^0(\spad,\Omm)=\{(f,F(f))\,|\, f\in \cS(\A_\K)_0/\overline{\{f-f_q\}}\}\oplus 2\C\oplus 2\C[1].
\end{equation}
\end{lem}
\proof
 Let   $ \alpha$ be a  nontrivial character of the additive group
$\A_\K/\K$. The lattice $\K\subset \A_\K$ is its own  dual. The  Fourier transform
\begin{equation}\label{fourier}
F(f)(a)=\int_{\A_\K} f(x)\alpha(ax)dx
\end{equation}
 depends on the choice of $\alpha$ but is canonical modulo the subspace $\{f-f_q\}$ and a fortiori modulo its closure $\overline{\{f-f_q\}}$.
 We recall that the Poisson formula  gives the equality
 $$
\sum_{q\in \K}f(q)=\sum_{q\in \K}F(f)(q)\qquad\forall f\in \cS(\A_\K)
$$
where $F$ denotes the Fourier transform. When applied to elements of $\cS(\A_\K)_0$ it gives
\begin{equation}\label{poisson2}
\sum_{q\in\K^\times} h(g^{-1}q)=|g|\sum_{q\in\K^\times} F(h)(gq).
\end{equation}
In particular, for $(f,h)\in \Ker\, \partial$ one gets
\begin{equation}\label{poisson3}
\Sigma(Fh)(g)=|g|^{-1}\Sigma(h)(g^{-1}),
\end{equation}
thus one obtains
\begin{equation}\label{coinv}
\Sigma(f-Fh)=0
\end{equation}
which shows that $f-Fh\in\overline{\{f-f_q\}}$,  by applying Lemma 5.4 of \cite{Meyer}.
\endproof\vspace{.05in}

\subsection{Spectral realization on   $H^1(\spad,\Omm)$}\label{rep}\hfill\vspace{.05in}

The action of the id\`ele class group $C_\K$ on the sheaf $\Omm$  is defined as follows
\begin{eqnarray}\label{rep1}
  \rep_+(\lambda)f(x)  &=& f(\lambda^{-1}x)\qqq f\in \Gamma(U_+,\Omm)\\
  \rep_-(\lambda)f(x)  &=& |\lambda|f(\lambda x)\qqq f\in \Gamma(U_-,\Omm) \nonumber \\
   \rep(\lambda)f(x)  &=& f(\lambda^{-1}x)\qqq f\in \Gamma(U_+\cap U_-,\Omm).  \nonumber
\end{eqnarray}\vspace{.05in}
To check that the equalities \eqref{rep1} determine a well-defined action of $C_\K$ on the sheaf $\Omm$, we need to show that \eqref{rep1} is compatible with the restriction maps \eqref{restmaps}. This is clear for the restriction from $U_+$ to $U$. For the restriction from $U_-$ to $U$, one has
\[
\Res(\rep_-(\lambda)f)(g)=|g|^{-1}\sum_{q\in \K^\times}(\rep_-(\lambda)f)(q g^{-1})=|g|^{-1}\sum_{q\in \K^\times}|\lambda|f(\lambda q g^{-1})
\]
while
$$
\rep(\lambda)\Res(f)(g)=\Res(f)(\lambda^{-1} g)=|\lambda^{-1} g|^{-1}\sum_{q\in \K^\times}f( q (\lambda^{-1} g)^{-1})
$$
and this shows the required compatibility.\endproof

 Let  $w$ be the element of the Weyl group $W$ of ${\rm PGL}_2$ given by the matrix
$$
w=
\left(
  \begin{array}{cc}
    0 & 1 \\
    1 & 0 \\
  \end{array}
\right).
$$
 This element acts on $C_\K$ by the automorphism $g\mapsto g^{-1}$. This determines the semi-direct product $N=C_\K\rtimes W$.  The element $w$ acts on $\spad$ by exchanging
 $0$ and $\infty$.

 We lift the action of $w$ to the sheaf $\Omm$ as follows. We consider the sheaf $w^{-1}_*\Omm$ which is defined by
\begin{equation}\label{winv}
    \Gamma(V,w^{-1}_*\Omm)=\Gamma(w(V),\Omm) \qquad\forall V \ {\rm open}.
\end{equation}
Then, we define the following morphism of sheaves $w_{\#}:\Omm\to w^{-1}_*\Omm$
\begin{eqnarray}\label{liftw}
  w_{\#} f  &=& f\in \Gamma(U_-,\Omm)\qqq f\in \Gamma(U_+,\Omm)\nonumber\\
 w_{\#} f  &=& f\in \Gamma(U_+,\Omm)\qqq f\in \Gamma(U_-,\Omm)  \\
   w_{\#} f(g)  &=& |g|^{-1} f(g^{-1})\qqq f\in \Gamma(U_+\cap U_-,\Omm).  \nonumber
\end{eqnarray}

The geometric action defined in the next proposition immediately implies the functional equation and in fact lifts the equation at the level of the representation

\begin{prop} \label{functequ}
1) The equalities \eqref{liftw} define an action of the Weyl group $W$ on the sheaf $\Omm$ which
fulfills the following compatibility property with respect to the action \eqref{rep1}
\begin{equation}\label{compinv}
    \rep(\lambda)w_{\#} \xi=|\lambda|w_{\#}\rep(\lambda^{-1})\xi.
\end{equation}
2) There is a unique action of the semi-direct product $N=C_\K\rtimes W$ on the sheaf $\Omm$ which agrees with \eqref{liftw} on $W$ and restricts on $C_\K$ to the twist $\rep[-\frac 12]$.
\end{prop}

\proof We need to show that the map $w_{\#}:\Omm\to w^{-1}_*\Omm$ defined in \eqref{liftw} is compatible with the restriction maps \eqref{restmaps}. For $f\in \Gamma(U_+,\Omm)$, one has
$$
\Res(w_{\#} f)(g)=|g|^{-1}\sum_{q\in \K^\times}f(q g^{-1})
$$
which agrees with $\displaystyle{w_{\#}\Res(f)}$, using \eqref{liftw}. A similar result holds for
$f\in \Gamma(U_-,\Omm)$. In fact, the full statement follows from the involutive property of the transformation $ w_{\#} f(g) = |g|^{-1} f(g^{-1})$.
\endproof

\begin{thm}\label{lemreasspec}
The cohomology $\displaystyle{H^1(\spad,\Omm)}$ gives the spectral realization of zeros of $L$-functions.
The spectrum of the action $\rep[-\frac 12]$ of $C_\K$ on $H^1(\spad,\Omm)$
 is invariant under the symmetry $\chi(g)\mapsto \chi(g^{-1})$ of Gr\"ossencharakters of $\K$.
\end{thm}

\proof
Consider the affine open set $U_-\subset \P^1_{\F_1}$
 and its complement $Y=\{0\}$. One checks directly that the cohomology with support $H_Y^1(\P^1_{\F_1},\Omm)$ gives the cokernel of the map $\Sigma:\cS_0(M)\to \strong(C_\K)$ of \eqref{Smap}, \ie the spectral realization of \cite{Meyer}, \cite{CCM},  \cite{CMbook}
 initiated in \cite{Co-zeta}.
 Since $U_-$ is affine, one can use the exact sequence \eqref{longexact1} which reduces to the isomorphism
\begin{equation}\label{court}
0\to H_Y^1(\P^1_{\F_1},\Omm)\to H^1(\P^1_{\F_1},\Omm)\to 0
\end{equation}
since the homomorphism
\begin{equation}\label{surject}
H^0(\P^1_{\F_1},\Omm)\to H^0(U_-,\Omm|_{U_-})
\end{equation}
is surjective \cf Lemma \ref{lemsurject}.
The symmetry then follows from the existence of the action of $N=C_\K\rtimes W$ on the sheaf $\Omm$ and hence on the cohomology
$\displaystyle{H^1(\spad,\Omm)}$.\endproof


\begin{thebibliography}{99}


\bibitem{Co-zeta} A.~Connes, {\em Trace formula in noncommutative geometry and the zeros of the Riemann zeta function}.  Selecta Math. (N.S.)  5  (1999),  no. 1, 29--106.

 \bibitem{ak} A.~Connes, C.~Consani  {\em On the notion of geometry over $\F_1$}, to appear in Journal of Algebraic Geometry; arXiv08092926v2 [mathAG].

\bibitem{CCM} A.~Connes, C.~Consani, M.~Marcolli, {\em
Noncommutative geometry and motives: the thermodynamics of endomotives},
Advances in Math. 214 (2) (2007), 761--831.

\bibitem{CCM2} A.~Connes, C.~Consani, M.~Marcolli, {\em The Weil proof
and the geometry of the adeles class space}, to appear in ``Algebra, Arithmetic
and Geometry -- Manin Festschrift'', Progress in Mathematics, Birkh\"auser
(2008); arXiv0703392.

\bibitem{ccm}  A.~Connes, C.~Consani, M.~Marcolli, {\em
Fun with $\F_1$}, Journal of Number Theory 129 (2009) 1532--1561.

\bibitem{CMbook} A.~Connes, M.~Marcolli,  {\em  Noncommutative Geometry, Quantum Fields, and Motives},
Colloquium Publications, Vol.55, American Mathematical Society, 2008.



 \bibitem{deit} A.~Deitmar {\em Schemes over F1}, in Number Fields and Function Fields –
Two Parallel Worlds. Ed. by G. van der Geer, B. Moonen, R. Schoof. Progr. in
Math, vol. 239, 2005.

 \bibitem{deit1} A.~Deitmar {\em  F1-schemes and toric varieties},  Contributions to Algebra and Geometry Vol. 49, No. 2, pp. 517-525 (2008).

\bibitem{deit2} A.~Deitmar {\em Remarks on zeta functions and K-theory over F1}. Proc. Japan Acad. Ser. A Math. Sci. 82 no 8, 141-146 (2006).

 \bibitem{demgab} M.~Demazure, P.~Gabriel {\em  Groupes alg\'ebriques}, Masson \& CIE, \'Editeur Paris 1970.

 \bibitem{Gilmer} R.~Gilmer, {\em Commutative semigroup rings} University of Chicago Press (1980).

 \bibitem{tohoku} A.~Grothendieck {\em Sur quelques points d'algèbre homologique} Tohoku Math. J., t. 9, 1957, p. 119- 183.





\bibitem{Hart} R.~Hartshorne  {\em Algebraic Geometry}, Graduate Texts in Mathematics 52, Springer-Verlag, New York Heidelberg Berlin 1977.

\bibitem{Ingham} A.~Ingham, {\em The distribution of prime numbers} With a foreword by R. C. Vaughan. Cambridge Mathematical Library. Cambridge University Press, Cambridge, 1990. xx+114 pp.

\bibitem{Kapranov} M.~Kapranov and A.~Smirnov, {\em Cohomology
determinants and reciprocity laws} Prepublication.

\bibitem{Kato} K.~Kato {\em Toric Singularities}
 American Journal of Mathematics, Vol. 116, No. 5 (Oct., 1994), 1073-1099.

\bibitem{Ku} N.~Kurokawa {\em Multiple zeta functions: an example.}
 in Zeta functions in geometry (Tokyo, 1990), Adv. Stud. Pure Math., 21, Kinokuniya, Tokyo, 1992, pp.  219--226.

\bibitem{KOW} N.~Kurokawa,  H.~Ochiai, A.~Wakayama, {\em Absolute Derivations
and Zeta Functions} Documenta Math. Extra Volume: Kazuya Kato's
Fiftieth Birthday (2003) 565-584.

\bibitem{Manin} Y.~I.~Manin, {\em Lectures on zeta functions and motives (according to Deninger and Kurokawa)} Columbia University Number-Theory Seminar (1992),
Ast\'erisque No.~228 (1995), 4, 121--163.

\bibitem{Meyer} R.~Meyer, {\em On a representation of the idele class group
related to primes and zeros of $L$-functions}.  Duke Math. J. Vol.127 (2005),
N.3, 519--595.

\bibitem{Soule} C.~Soul\'e, {\em Les vari\'et\'es sur le corps \`a un
\'el\'ement}. Mosc. Math. J. 4 (2004), no. 1, 217--244.

\bibitem{Steinberg} R.~Steinberg {\em A geometric approach to the representations of the full linear group over a Galois field} Transactions of the AMS, Vol. 71, No. 2 (1951), pp. 274--282.



\bibitem{Tits} J.~Tits, {\em Sur les analogues alg\'ebriques des groupes semi-simples
complexes}. Colloque d'alg\`ebre sup\'erieure, Bruxelles 19--22
d\'ecembre 1956, Centre Belge de Recherches Math\'ematiques \'Etablissements
Ceuterick, Louvain; Librairie Gauthier-Villars, Paris   (1957), 261--289.



\bibitem{TV} B.~T\"oen, M.~Vaqui\'e {\em Au dessous de $\text{Spec}(\Z)$}
 (preprint) arXiv0509684v4 to appear in K-theory.



\end{thebibliography}
\end{document}